\newif\ifprint
\printfalse
\documentclass[twoside,10pt]{HandbookOfModuli}

\usepackage{amsmath,amsthm}
\usepackage{amssymb}
\usepackage{amscd}
\usepackage{latexsym}
\usepackage[utf8]{inputenc}
\usepackage{textcomp}
\usepackage{color}
\usepackage[pdftex]{graphicx}

\usepackage{titlesec}
\usepackage{xspace}
\usepackage[all]{xypic}
\ifprint
	\input{giovanni} 
	\usepackage[final]{microtype}
\fi
\usepackage{bm}
\renewcommand{\mathbf}[1]{\bm{#1}} 
\ifprint
	\definecolor{linkred}{rgb}{0,0,0} 
	\definecolor{linkblue}{rgb}{0,0,0} 
\else
	\definecolor{linkred}{rgb}{0.7,0.2,0.2}
	\definecolor{linkblue}{rgb}{0,0.2,0.6}
\fi

\numberwithin{equation}{section} 

\usepackage[
	hypertexnames=false,
	hyperindex,
	pagebackref,
	pdftex,
	pdftitle={Handbook of Moduli},
	pdfdisplaydoctitle,
	pdfpagemode=UseNone,
	breaklinks=true,
	extension=pdf,
	bookmarks=false,
	plainpages=false,
	colorlinks,
	linkcolor=linkblue,
	citecolor=linkred,
	urlcolor=linkred,
	pdfmenubar=true,
	pdftoolbar=true,
	pdfpagelabels,
	pdfpagelayout=TwoPage,
	pdfview=Fit,
	pdfstartview=Fit
]{hyperref}

\catcode`@=11 \baselineskip 4.5mm
\def\ps@handbook{\def\@oddhead{\hfill \leftmark \hfill\thepage }
\def\@evenhead{\thepage \hfill \rightmark \hfill}
\def\@oddfoot{}
\def\@evenfoot{}}

\def\@evenhead{}
\def\@oddfoot{}
\def\@evenfoot{\hfill\copyright\ China Higher Education Press}
\def\list#1#2{\ifnum \@listdepth >5\relax \@toodeep \else \global
\advance \@listdepth\@ne \fi \rightmargin \z@ \listparindent\z@
\itemindent\z@ \csname @list\romannumeral\the\@listdepth\endcsname
\def\@itemlabel{#1}\let\makelabel\@mklab \@nmbrlistfalse #2\relax
\@trivlist \parskip -\parsep \parindent\listparindent \advance
\linewidth -\rightmargin \advance\linewidth -\leftmargin \advance
\@totalleftmargin \leftmargin \parshape \@ne \@totalleftmargin
\linewidth \ignorespaces}

\renewcommand*\l@section{\@tocline{1}{0pt}{0em}{1em}{}}\renewcommand*\l@subsection{\@tocline{2}{0pt}{1.5em}{2em}{}} \renewcommand\contentsnamefont{\bfseries\large} \catcode`@=12
\renewcommand{\theequation}{\thesection.\arabic{equation}}

\def\thebibliography#1{\section*{References}
\list{[\arabic{enumi}]}{\settowidth \labelwidth{[#1]} \leftmargin
\labelwidth \advance \leftmargin \labelsep \usecounter{enumi}}
\def\newblock{\hskip .11em plus .33em minus .07em} \sloppy
\clubpenalty 4000 \widowpenalty 4000 \sfcode`\.=1000 \relax}

\titleformat{\section}{\normalfont\large\bfseries}{\thesection.}{0.5em}{}[\kern0.em]
\titleformat{\subsection}{\normalfont\bfseries}{\thesubsection.}{0.3em}{}[\kern0.em]
\titleformat{\subsubsection}[runin]{\normalfont\bfseries}{\thesubsubsection.}{0.5em}{}[\kern0.5em]

\def\fofsubsubsection#1{\refstepcounter{equation}\subsubsection*{\theequation.\kern0.25em #1}}
\def\foisubsubsection#1{\refstepcounter{equation}\subsubsection*{\kern\parindent\theequation.\kern0.25em #1}}

\newcommand{\gquot}{/\!\!/}
\newcommand{\Gm}{\mathbb{G}_m}

\newcommand{\bA}{\mathbb{A}}
\newcommand{\bP}{\mathbb{P}}
\newcommand{\bZ}{\mathbb{Z}}
\newcommand{\bQ}{\mathbb{Q}}
\newcommand{\bR}{\mathbb{R}}
\newcommand{\bC}{\mathbb{C}}
\newcommand{\calB}{\mathcal{B}}
\newcommand{\calH}{\mathcal{H}}
\newcommand{\calD}{\mathcal{D}}
\newcommand{\calh}{\mathfrak{h}}
\newcommand{\calM}{\mathcal{M}}
\newcommand{\calN}{\mathcal{N}}
\newcommand{\calT}{\mathcal{T}}
\newcommand{\calO}{\mathcal{O}}
\newcommand{\calL}{\mathcal{L}}
\newcommand{\calI}{\mathcal{I}}
\newcommand{\calP}{\mathcal{P}}

\DeclareMathOperator{\st}{st}
\DeclareMathOperator{\Spec}{Spec}
\DeclareMathOperator{\Proj}{Proj}
\DeclareMathOperator{\Pic}{Pic}
\DeclareMathOperator{\NS}{NS}
\DeclareMathOperator{\Hom}{Hom}

\DeclareMathOperator{\Nef}{Nef}
\DeclareMathOperator{\Mov}{Mov}
\DeclareMathOperator{\Cox}{Cox}
\DeclareMathOperator{\Sym}{Sym}
\DeclareMathOperator{\Hilb}{Hilb}
\DeclareMathOperator{\PGL}{PGL}
\DeclareMathOperator{\GL}{GL}
\DeclareMathOperator{\SL}{SL}

\numberwithin{equation}{section}
\theoremstyle{plain}
\newtheorem{theorem}[equation]{Theorem}
\newtheorem{lemma}[equation]{Lemma}
\newtheorem{proposition}[equation]{Proposition}

\theoremstyle{remark}
\newtheorem{remark}[equation]{Remark}

\theoremstyle{definition}
\newtheorem{definition}[equation]{Definition}
\newtheorem{notation}[equation]{Notation}

\begin{document}
\include{HandbookOfModuliPreamble}
\bibliographystyle{amsalpha}

\title[GIT and moduli]{GIT and moduli with a twist}
\author[R. Laza]{Radu  Laza}
\address{Stony Brook University, Department of Mathematics,  Stony Brook, NY 11794}
\email{rlaza@math.sunysb.edu}
\thanks{The author is partially supported by NSF grant DMS-0968968 and a Sloan fellowship.}


\begin{abstract}
We survey the role played by GIT in the study of moduli spaces, with an emphasis on the birational geometry of GIT quotients.
\end{abstract}
\maketitle

\section{Moduli and GIT -- a brief history} 
Geometric Invariant Theory  (GIT) is an important tool in the study of moduli spaces in algebraic geometry. In fact,  Mumford said in the preface of the first edition of the foundational GIT book \cite{GIT}: ``{\it to construct moduli schemes for various types of algebraic objects ... appears to be, in essence, a special and highly non-trivial case}'' of  GIT. The modern point of view is more nuanced and in many instances moved away from GIT.  For instance, stacks, already used in  \cite{dm},  have become the language of modern moduli theory. Furthermore, their use is essential in many situations (e.g. \cite{av} from which we inspired our title).  In another direction,  the ideas pioneered by Koll\'ar, Shepherd-Barron,  and Alexeev (\cite{ksb}, \cite{alexeev}) and based on the minimal model program give an approach to the compactification of moduli spaces of higher dimensional varieties without using GIT (see the survey \cite{kollar} in this volume).  Finally, the  log structures enlarge the notion of smoothness so that some moduli spaces are naturally compact (see the survey \cite{abramovich} in this volume).  Yet, arguably GIT still plays an important role in moduli theory today. The twist in the title alludes to the developments since the appearance of variation of GIT quotient theory (abbreviated as VGIT in what follows), due to Dolgachev--Hu \cite{dh} and Thaddeus \cite{thaddeus}, in which  the birational geometry of GIT quotients (and moduli spaces) plays an increasingly central role.  This paper is a modest attempt to survey these developments.

The first instance of GIT and moduli spaces is probably the study of the moduli space of elliptic curves. Since an elliptic curve can be embedded as a cubic curve in $\bP^2$ uniquely up to projective equivalence, an algebraic description of the moduli space is the (GIT) quotient $\bP\Sym^3 V^*\gquot \PGL(3)$  (with $V$  the standard $\PGL(3)$ representation). It turns out that $\bP\Sym^3 V^*\gquot\PGL(3)\cong \bP^1$ corresponding to the fact that the ring of $\PGL(3)$-invariant polynomials is isomorphic to the polynomial algebra $k[S,T]$ (with the standard $j$-invariant being the rational function $\frac{16S^3}{T^2+64S^3}$). The investigation  of this type of problems, i.e. finding the invariants for various group actions coming from geometry, was one of the most active areas in mathematics in the nineteenth century.  The highlights of that period include the computations of  invariants for $n\le 6$ points on $\bP^1$, quartic curves, and cubic surfaces. Finding explicit invariants is a difficult task, and our knowledge today is not much better (see however \cite{vakil} for a discussion of the ring of invariants for $n$ ordered points on $\bP^1$). This classical period of invariant theory essentially ended with the arrival (around 1900) of the fundamental result of Hilbert that says that the ring $\Sym(V^*)^{\PGL(n)}$ of invariant polynomials for a linear representation $V$ is finitely generated. 

The subject of invariant theory was reinvented by Mumford in the sixties. The work of Mumford  \cite{GIT} put the subject on firm theoretical footing and showed how to use GIT without explicitly knowing the invariants. In particular, Mumford used GIT to show that the moduli space of curves $M_g$ is quasi-projective.  Later, Mumford and Gieseker (\cite{mumford}) proved that  (the coarse moduli space associated to) the Deligne--Mumford compactification $\overline{M}_g$ of the moduli space of genus $g$ curves is a projective compactification of $M_g$ that can be constructed via GIT. Some other major results around the same time include the proof of quasi-projectivity for the moduli of surfaces of general type (Gieseker \cite{gieseker}) and compactifications for the moduli spaces of vector bundles over curves (Mumford, Narasimhan, Seshadri, e.g. \cite{seshadri}) and surfaces (Gieseker \cite{giesekervb}), and then torsion free sheaves (Maruyama \cite{maruyama}). 

After a very active period in the sixties and seventies, the focus in moduli theory and GIT shifted somewhat. Among the more important later results we mention the work of Kirwan on the cohomology of certain moduli spaces via GIT (\cite{kirwancoh}) and on partial desingularizations of GIT quotients (\cite{kirwan}). Also, Viehweg \cite{viehweg} proved that the moduli space of varieties of general type exists as a quasi-projective variety. It is important to note that  the work of Viehweg  does not continue the approach of Mumford and Gieseker for curves and surfaces. Namely, in order to prove stability for generic smooth varieties, Viehweg used a non-standard polarization on the Hilbert scheme, in effect prescribing the set of stable points. The changing of the linearization in a GIT problem is nowadays a common occurrence, but it has become so only in the next phase of the GIT story: the theory of variation of GIT quotients, discussed in the following paragraph.

A renewal of GIT occurred in the early nineties starting from a well-known, but essentially ignored, issue in the construction of GIT: the construction of a quotient $X\gquot G$ depends on the choice of a $G$-linearized ample line bundle $\calL$ on $X$. Thaddeus \cite{thaddeus}
 and independently Dolgachev--Hu \cite{dh} (and previously Brion--Procesi \cite{brion} in the toric case) analyzed this dependence and showed that it is surprisingly well behaved.   Roughly speaking, there are only finitely many distinct possibilities for the GIT quotients, which are related by quite explicit birational operations.  A little later, Hu--Keel \cite{hukeel}   showed that these results are connected in a fundamental way with the most well behaved spaces in the minimal model program, the so-called  {\it Mori dream spaces}.
 
 The influence of the ideas introduced by the arrival of VGIT extends beyond the direct applications of the theory to moduli spaces. Namely,  the surprisingly nice results of VGIT together with substantial progress in  birational geometry (see \cite{bchm}), suggest that the various birational models of a moduli space are related in a meaningful and useful way. More specifically, different constructions (or different choices in a construction) for a moduli space typically give different compactifications. Initially, this was thought of as a pathological feature of moduli theory. VGIT contributed to a significant shift in this view: one realized  that the existence of multiple compactifications can be used to one's advantage. Typically, using a compactification with relatively simple structure, one can extract important information (e.g. cohomology) about other compactifications of more geometric interest.  In other words, each compactification gives a facet of the moduli problem and taken together one gets a fuller picture. There are numerous instances of this principle in the current literature  on moduli spaces, but we choose to mention here only  two examples familiar to the author: (1) the search for a canonical model  for the moduli space of curves $\overline{M}_g$ (e.g. \cite{hh,hh2}, and the survey \cite{fs} in this volume), and (2) the comparison between certain period domains  and GIT quotients (e.g. \cite{lc1,lc2}). 

\subsubsection*{Content.} After a brief discussion (Sect. \ref{standardgit}) of the constructions and results in the standard GIT situation, we survey (Sect. \ref{secresults}) the main results of the theory of variation of GIT quotients. We then (Sect. \ref{sectools}) discuss some tools (such as the numerical criterion) that can be used to understand GIT quotients in general, but focus on the VGIT situation. In Section \ref{sectbir}, we review the relationship between VGIT and birational geometry.  We close with a survey of some applications of the birational geometry of GIT quotients and moduli spaces (Sect. \ref{secapp}).

\subsubsection*{Acknowledgement.} The author is grateful to   M. Thaddeus and I. Dolgachev from whom he learned about the subject.  The comments of M. Fedorchuk, N. Giansiracusa, K. Schwede, and D. Swinarski have helped improve the manuscript. 

\subsubsection*{Disclaimer.} The omissions and inaccuracies are solely the responsibility of the author. Also, the topics included are not intended to  exhaust the subject, but rather reflect the interests and expertise of the author. In particular, one important topic not discussed here is the connection between GIT/VGIT and the moment map (see \cite[Ch. 8]{GIT} and \cite{dh}).

\subsubsection*{Conventions.} We work over an algebraically closed field $k$ of characteristic $0$.  Throughout the paper $G$ is a reductive group acting on a  quasi-projective variety $X$. The variety is always assumed to be normal, and quite frequently smooth. 
We will denote orbits by $G\cdot x$ and stabilizers by $G_x$. $G^0$ (and $G_x^0$) will denote the connected component of the corresponding group.

\section{The GIT construction and main results}\label{standardgit}
This section aims to give a brief overview of the standard GIT as developed by Mumford \cite{GIT}. Other standard textbook references are \cite{newstead},  \cite{popov}, \cite{dolgachev}, and \cite{mukaib}. 
\subsection{Affine Quotients.}\label{sectaff}
The simplest instance of a GIT quotient is that of a reductive group $G$ acting linearly on an affine variety $X=\Spec R$. In this situation, a foundational result of Hilbert says that the ring of invariants $R^G\subset R$ is a finitely generated $k$-algebra (\cite[Thm. A.1.0]{GIT}, \cite[Thm. 3.2, 3.3]{dolgachev}). Thus, it is natural to define the quotient to be the affine variety:
\begin{equation}\label{defaff}
X/G:=\Spec R^G.
\end{equation}

\begin{remark}
The assumption of reductive group is essential for the finite generation of $R^G$ (cf. Nagata \cite{nagata}; see \cite{mukai} for some geometric counter-examples, and \cite{knonred} for a survey of results on quotients by non-reductive groups). Two special cases of reductive groups are of particular interest: $G$ is a semi-simple group (e.g. $\PGL(n)$) or a  torus $(\mathbb G_m)^n$. In general a reductive group is built out of these two cases  (see \cite[App. A]{GIT} for more on reductive groups).  
\end{remark}

\subsubsection{Categorical vs. Geometric Quotient.} In the absence of additional choices (such as characters for $G$, as discussed elsewhere in this survey), the definition \eqref{defaff} is essentially the only possibility in the realm of algebraic geometry. More precisely, the natural $G$-invariant projection $\pi:X\to X/G$ makes $X/G$ a {\it universal categorical quotient} (see \cite[Def. 0.7]{GIT} and \cite[Thm. 1.1]{GIT}), i.e. any $G$-invariant morphism $f:X\to Y$ factors through $X/G$:
$$
\xymatrix
{ 
X\ar@{->}[r]^{\pi}\ar@{->}[rd]_{f}&X/G\ar@{->}[d]^{\exists! \bar f}\\
&Y\\
}.
$$

However,   $X/G$  is typically not a {\it geometric quotient} (see \cite[Def. 0.6]{GIT}). In general, there is no one-to-one correspondence between the points of $X/G$ and the orbits of $G$. For a simple example, consider the action of $G=\bC^*$ on $X=\bA_{\bC}^1\cong \bC$ given by the natural multiplication $(t,x)\in \bC^*\times \bC\to t\cdot x\in \bC$. Since $\bC[x]^G\cong \bC$ (the only invariants are the constants), the quotient $X/G$ is a point and the quotient map $\pi$ is the trivial projection $\bC\to \{\mathrm{pt.}\}$. On the other hand, the action of $G$ on $X$ has two orbits: $\{0\}=G\cdot\{0\}$ and $\bC\setminus \{0\}=G\cdot \{1\}$. The issue here is that the fibers of $\pi$ are always closed in $X$. In our example, the orbit $\{0\}$ is closed, while the closure of the orbit $\bC\setminus \{0\}$ is the affine line, which contains the closed orbit $\{0\}$. Thus, the two orbits map to the same point via $\pi$, showing  that the orbits are not always separated by the quotient map.

\subsubsection{Closed Orbits.}\label{closedorb} In some sense, when $G$ is reductive, the failure of $X/G$ to be a geometric quotient is always of the type exemplified above. Namely, we recall the following fact about the orbits of algebraic group actions (e.g. \cite[\S8.3]{hum}):  {\it each $G$-orbit is smooth, locally closed in $X$, and its boundary is a union of orbits of strictly lower dimension}. This easily gives that  each fiber of $\pi$ contains a unique closed orbit, namely, the orbit of minimal dimension in that fiber. Furthermore, if $G\cdot x_0$ is a closed orbit in $X$, then for all $x\in \pi^{-1}(\pi(x_0))$, the orbit of $G\cdot x$ contains $G\cdot x_0$ in its closure. 

\subsection{Projective Quotients.}\label{sectproj} 
One is typically  interested  in the action of a  reductive group $G$ on a projective variety $X$. Mumford constructed a GIT quotient in this situation by  considering an ample line bundle $\calL$ together with a $G$-linearization (i.e. essentially a lift of the $G$-action from $X$ to $\calL$; see \cite[\S3]{GIT}). This choice gives an  embedding: 
$$i:X\xrightarrow{|\calL^{\otimes k}|} \bP^N \ \ (\textrm{for some $k\gg0$})$$ such that $G$ acts linearly on $\bP^N$ and the embedding $i$ is $G$-equivariant. By considering affine cones, one can reduce to the affine situation. Concretely, the definition of a GIT quotient\footnote{Alternatively, the general GIT quotient can be defined by gluing affine quotients (see the proof of \cite[Thm. 1.10]{GIT}). In the situation discussed here (esp. $\calL$ is ample), the two approaches are equivalent.} in the projective case is as follows.

\begin{definition}
Let $G$ be a reductive group acting on a projective variety $X$. For $\calL$ an ample $G$-linearized line bundle, the associated {\it GIT quotient} is the projective variety:
\begin{equation}\label{defgit}
X\gquot_\calL G:=\Proj \oplus_{n\ge 0} H^0(X,\mathcal L^{\otimes n})^G.
\end{equation}
\end{definition}
As in the affine case, one shows that $X\gquot G$ has good categorical properties and thus the definition \eqref{defgit} is in some sense canonical (see \cite[Thm. 1.10]{GIT}). 
 
 \begin{remark}\label{remrelative}
 A more general situation is the relative situation: $X$ projective over an affine variety $Y$, $G$ acts equivariantly on $X\to Y$, and $\calL$ is a relative ample line bundle on $X$. The above definitions can be easily adapted to this relative situation. In particular, we note that there exists a structural morphism $X\gquot_ \calL G\to Y/G$, where $X\gquot _\calL G$ is a projective quotient (the relative version of \eqref{defgit}) and $Y/G$ is an affine quotient (as in \eqref{defaff}). A particularly interesting case is the toric case discussed in \S\ref{sectoric}: $X=Y=\bA^n$, $G=T=(\Gm)^r$, $\calL=\calO_X$.
 \end{remark}
 
\subsubsection{Semistable points.} A new type of behavior appears in the case of projective quotients. Namely, 
the natural quotient map $\pi:X\to X\gquot_{\calL} G$ is only a rational map:  the domain of definition of $\pi$ is precisely the set of points $x\in X$ such that there exists a $G$-invariant section  $\sigma\in H^0(X,\mathcal L^{\otimes n})^G$ (for some $n$) with $\sigma(x)\neq 0$ (N.B. here $X_\sigma$ is automatically affine). We call such points {\it semistable} and denote by $X^{ss}(\calL)\subset X$ the corresponding open set. The points in $X^{us}(\calL):=X\setminus X^{ss}(\calL)$ are called {\it unstable} and are excluded from the GIT analysis. 

\subsubsection{Stable points.} 
In moduli theory, one is particularly interested in geometric quotients. The discussion of \S\ref{closedorb} applies also here.  It follows that for points $x\in X^{ss}(\calL)$ such the orbit $G\cdot x$ is closed in $X^{ss}$ and of maximal dimension (i.e. $\dim G\cdot x=\dim G$, or equivalently the stabilizer $G_x$ is finite) one has $\pi^{-1}(\pi(x))=G\cdot x$. One calls such points {\it stable points}. 
 The set of stable points $X^s(\calL)\subset X^{ss}(\calL)$ is an open $G$-invariant subset such that the induced quotient $X^s(\calL)/G$ is both a geometric and categorical quotient. If $X^s(\calL)\neq \emptyset$, then $X^s(\calL)/G$ is an open dense subset of the quotient $X\gquot_\calL G$. 
 
\begin{remark} To emphasize that the quotient map $X\dashrightarrow X\gquot G$ is only defined on $X^{ss}$ and that $X\gquot G$ is only a categorical quotient, one uses the notation $\gquot$. In contrast, $X^s\to X^s/G$ satisfies all expected properties and thus the  notation $/$.
\end{remark}

\subsection{GIT and moduli} GIT and moduli theory are closely related since in many situations it is possible to construct good parameter spaces $X$ (e.g. Hilbert  and Quot schemes) for algebraic objects (with extra rigidifying structure) on which a reductive group $G$ acts naturally. For instance, one might be interested in (smooth) algebraic varieties of a certain type T. A typical first step in constructing a moduli space for them is to prove that all (smooth) varieties of type T have a uniform embedding  $V\hookrightarrow \bP^N$ in a large projective space (e.g. smooth curves embedded by some large power $\omega_C^{\otimes \nu}$ of the canonical bundle). If this is true, then the smooth varieties of type $T$ are parameterized by a locally closed subset of some irreducible component $X$ of $\Hilb_p(\bP^N)$, where $p$ is the corresponding Hilbert polynomial of $V$. The group $G=\PGL(N+1)$ acts naturally on $X$ by change of coordinates on $\bP^N$ (which amounts to changing the embedding $V\hookrightarrow \bP^N$). Thus, a moduli space for varieties of type T would be roughly the GIT quotient $X\gquot G$, which is a projective variety. Unfortunately,  it is very hard to prove that even the generic variety of type T is GIT stable. Successful applications of GIT to constructions of moduli spaces with quasi-projective coarse scheme include the important cases of abelian varieties (\cite[Thm. 7.10]{GIT}) and curves (\cite[Cor. 7.14]{GIT}). For curves, it is possible to control the stability enough to obtain that the Deligne--Mumford compactification has a projective coarse moduli space (\cite[Thm. 5.1]{mumford}; see \cite{msurvey} for a survey).

The most desirable case is when $X^s(\calL)=X^{ss}(\calL)$. Namely, one gets that the quotient $X\gquot_\calL G=X^s(\calL)/G$ is both a projective variety and a geometric quotient. In good situations (e.g. $\overline{M}_g$), this quotient has also a modular interpretation.   If $X$ is smooth, the quotient $X\gquot_\calL G$, considered as a stack,  is a smooth proper Deligne--Mumford stack with a projective coarse scheme (arguably, this is an ideal outcome for a moduli problem). As we will see later,  quite often, $X^s(\calL)=X^{ss}(\calL)$ happens for all generic choices of $\calL$. Unfortunately, the natural GIT set-up for many moduli problems gives  situations with $X^s(\calL)\subsetneq X^{ss}(\calL)$. Even in those situations the fact that $X\gquot_\calL G$ is projective can be used to one's advantage. Specifically, the properness of the quotient gives the following useful {\it semistable replacement property} (\cite[Lem. 5.3]{mumford}, \cite[Prop. 2.1]{shah}).
 \begin{lemma}\label{ssreplace}
Let $S=\Spec R$ and $S^*=\Spec (K)$, where $R$ is a DVR with  field of fractions $K$ and closed point $o$. Assume that $S^*\to X^s/G$ for some GIT quotient. Then, after a finite base change $S'\to S$ (ramified only at the special point $o$), there exists a lift $\tilde{f}:S'\to X^{ss}$ of $f$ as in the diagram:
$$
\xymatrix
{ 
&S'\ar@{->}[r]^{\tilde{f}}\ar@{->}[ld]_{f}&X^{ss}\ar@{->}[rd]\\
S&S^*\ar@{_{(}->}[l]\ar@{->}[r]&X^s/G\ar@{^{(}->}[r]&X\gquot G\\
}.
$$
 Furthermore, one can assume that $\tilde{f}(o)$ belongs to a closed orbit. 
\end{lemma} 

In other words, while a GIT quotient typically fails to have a modular meaning at the boundary, one can use this lemma  to understand the degenerations of smooth objects and then construct or understand a good compactification of the moduli space (see \cite{shah} and \cite{caporaso} for some concrete applications of this principle). For a formalization, from the perspective of stacks, of the properties satisfied by the GIT quotients occurring in constructions of moduli spaces see \cite{alpergit}.

\subsection{Some concluding remarks on standard GIT}
For  a reductive group $G$ acting on a projective variety $X$, the discussion above can be summarized as:
\begin{itemize}
\item[(1)] The construction of a quotient is based on the fact that the ring of invariants $R^G$ is finitely generated. The quotient has good functorial properties.
\item[(2)] $X\gquot G$ is a projective variety. 
\item[(3)] The quotient map $\pi$ is defined only on the semistable locus $X^{ss}$.
\item[(4)] Each fiber of  $\pi: X^{ss}\to X\gquot G$ contains a unique  closed orbit in $X^{ss}$. Furthermore, $\pi(x)=\pi(y)$ iff $\overline{G\cdot x}\cap \overline{G\cdot y}\cap X^{ss}\neq \emptyset$.
\item[(4)] $X^{s}/G$ is a geometric quotient; it is an open and, if non-empty, dense subset of $X\gquot G$. In particular, $G_x$ is finite and $\pi^{-1}(\pi(x))=G\cdot x$ for $x\in X^{s}$. 
\item[(5)] In the cases coming from moduli problems: if   $X^{ss}=X^s$, then  the corresponding moduli space (e.g. $\overline{M}_g$) is a  Deligne--Mumford stack (and smooth if $X$ is smooth) with the associated  coarse scheme being a projective normal variety . Even when $X^{ss}\subsetneq X^s$, the projectivity of the GIT quotient is useful to analyze the degenerations of smooth objects and to compactify the moduli space.
\end{itemize} 

While the definition  \eqref{defgit} of the GIT quotient depends on the choice of linearization $\calL$,  we choose to ignore $\calL$ from notation in this summary to emphasize that for a long time the dependence on $\calL$ was essentially ignored (see however \cite{gm}). One reason for this might be that in many GIT situations coming from moduli problems there is a preferred choice for $\calL$ (e.g. asymptotic linearizations on Hilbert schemes). A systematic investigation of the dependence on the linearization started with the toric case (\cite{gm}) and culminated with the  results of Dolgachev--Hu  \cite{dh} and Thaddeus  \cite{thaddeus}  discussed below.

\section{The main results of VGIT}\label{secresults}
In this section we review the main results of the theory of variation of GIT quotients following Thaddeus \cite{thaddeus} and Dolgachev-Hu \cite{dh}.

\subsection{The space of linearizations and the partition into chambers} 
The first step in understanding the dependence of the GIT quotient $X\gquot_\calL G$ on the linearization $\calL$ is to define a GIT equivalence relation for linearizations.  While the space of linearizations  is essentially already understood in \cite[\S1.3]{GIT}, the fact that the natural  GIT equivalence  is well behaved (cf. \cite[Thm. 2.4]{thaddeus}, \cite[Thm. 0.2.3]{dh}) is a surprising result and one of the cornerstones of VGIT.

\subsubsection{The parameter space for linearizations.} A linearization consists of an underlying line bundle together with the extra data of  a $G$-linearization (\cite[Def. 1.6]{GIT}). Thus, denoting by $\Pic^G(X)$ the space of $G$-linearized line bundles, there is a forgetful map $\Pic^G(X)\to \Pic(X)$, whose kernel is $\chi(G)$, the group of characters of $G$. However, the GIT quotient $X\gquot_\calL G$ and the sets of stable and semistable points $X^{s(s)}(\calL)$ only depend on the algebraic equivalence class of the linearization $\calL$ (see \cite[\S2]{thaddeus}, \cite[Def. 2.3.4]{dh} for definitions). More precisely, one has (cf. \cite[Prop. 2.1]{thaddeus}, \cite[Prop. 2.3.6]{dh}):
\begin{proposition}\label{proplinear}
If $\calL$ is an ample linearization, then $X^{ss}(L)$, and the quotient $X\gquot_\calL G$ regarded as a polarized variety, depend only on the $G$-algebraic equivalence class of $\calL$.
\end{proposition}

It follows that  the relevant parameter space for linearizations is a finitely generated abelian group  $\NS^G(X)$, called {\it the $G$-linearized Neron--Severi group}, which  has a natural forgetful map $\NS^G(X)\to \NS(X)$ to the usual Neron--Severi group.  Note however that not every line bundle $\calL$ can be linearized, but, if $X$ is normal\footnote{The situation in the non-normal case is significantly more delicate, e.g.  \cite[\S4]{alexeevab}.}, some power $\calL^{\otimes n}$ can be linearized (cf. \cite[Cor. 1.6]{GIT}). Also,  changing a $G$-linearized bundle $\calL$ by a multiple does not change $X^{ss}(\calL)$ or the quotient $X\gquot_\calL G$; it only changes the polarization of $X\gquot_\calL G$ by the same multiple. Thus, it is preferable to work with numerical equivalences and with $\bQ$ coefficients ({\it fractional linearizations} in \cite{thaddeus}), i.e. to consider $\NS^G_\bQ(X):=\NS^G(X)\otimes_\bZ \bQ$ as the parameter space for linearizations. One gets (see \cite[\S2.3.9]{dh} for $\bZ$ coefficients):

\begin{proposition}
Let $G$ be a reductive group acting on a normal projective variety $X$. Then, the space of linearizations $\NS_\bQ^G(X)$ sits in an exact sequence:
$$0\to \chi(G)\otimes_\bZ \bQ\xrightarrow{i} \NS_\bQ^G(X)\xrightarrow{f} \NS(X)\otimes_\bZ\bQ\to 0,$$
where the inclusion $i$ corresponds to the choices of linearization on the trivial bundle $\calO_X$ and $f$ is the map that forgets the data of $G$-linearization. In particular, $\NS_\bQ^G(X)$ is a finite dimensional vector space of dimension $\rho(X)+\tau(G)$, where $\rho(X)$ is the Picard number of $X$ and $\tau(G)$ is the dimension of the radical of $G$ (a torus).
\end{proposition}

\begin{remark}
For $G=\SL(n)$, every line bundle carries a unique $G$-linearization (N.B. $\PGL(n)$ can be understood via the isogeny to $\SL(n)$, see \cite[p. 33]{GIT}). Thus, the choice of linearization is equivalent to the choice of an ample line bundle on $X$. On the other hand, for a torus $G=(\Gm)^n$, one has a lattice of characters $M:=\chi(G)\cong \bZ^n$. A particular case of interest for the toric case is when $X=\bA^m$ is an affine variety. Thus, the only bundle is $\calO_X$, and $\NS^G(X)\equiv M$. Quotients of this type describe the toric varieties (see \S\ref{sectoric}).
\end{remark}

As is standard in algebraic geometry, one considers the convex cone $A(X)$ spanned by ample line bundles in $N^1_\bR(X)$  (where $N^1_\bR(X)=\NS(X)\otimes_\bZ \bR$ are the numerical equivalence classes with real coefficients). By definition, the GIT quotient  only makes sense on the cone $f^{-1}(A(X))\cap \NS^G_\bQ(X)$, where $f:\NS^G_\bR(X)\to N^1_\bR(X)$ is the natural forgetful map. Furthermore, there might exist ample $G$-linearized bundles $\calL$ for which  there are no semistable points (and thus the quotient is empty).  It follows that one has to restrict to the {\it cone of  $G$-ample line bundles}: 
\begin{equation}\label{eqcone}
C^G(X)\subseteq f^{-1}(A(X))\subset \NS_\bR^G(X)
\end{equation}
spanned by classes of ample $G$-linearized line bundles $\calL$ such that $X^{ss}(\calL)\neq \emptyset$ (i.e. $\calL$ is ample and $G$-effective, or simply $G$-ample, cf.  \cite[p. 701]{thaddeus}, \cite[Def. 3.1.1, Def. 3.2.1]{dh}). A first result of VGIT is then (\cite[(2.3)]{thaddeus}):

\begin{proposition}
The cone $C^G(X)$ is convex and it is the intersection of  $f^{-1}(A(X))$ with a rational polyhedron in $\NS_\bR^G(X)$. 
\end{proposition}

The boundary of the cone corresponds to either  ample line bundles for which $X^s(\calL)=\emptyset$ (and thus the quotient $X\gquot G$ has less than expected dimension) or nef but not ample line-bundles (\cite[Prop. 3.2.8]{dh}). In particular, note that if $X^s(\calL)\neq \emptyset$ for some ample $\calL$, the interior of the cone $C^G(X)$ is not empty. Note also that $C^G(X)$ is not necessarily rational polyhedral (it is so when regarded as a subset of $f^{-1}(A(X))$). This phenomenon has to do with the pathological behavior of the ample cone of $X$. For example, one can take $G=\{1\}$ (the trivial group) and $X$ one of the varieties with ``round'' ample cone (then $C^G(X)=A(X)$ is round, but the proposition is still valid). This shows that the structural results of VGIT have to do more with $G$ than $X$.  

\subsubsection{GIT equivalence and the partition of $C^G(X)$ into chambers.} We have discussed above that the parameter space for linearizations is $C^G(X)\cap \NS^G_\bQ(X)$. Then, one defines a coarsening of the  algebraic equivalence of linearizations (compare Prop. \ref{proplinear}) as follows (\cite[Def. 3.4.1]{dh}).  

\begin{definition}
For $G$ and $X$ as before, we say that the linearizations $\calL$ and $\calL'$ are {\it GIT equivalent} if $X^{ss}(\calL)=X^{ss}(\calL')$ (which also implies that $X^s(\calL)=X^s(\calL')$ and $X\gquot_\calL G\cong X\gquot_{\calL'} G$). 
\end{definition}

A first result of VGIT  is that there are only finitely many GIT equivalence classes, i.e. there are only finitely many open subsets of $X$ that can be realized as $X^{ss}(\calL)$ for some linearization (cf. \cite[Thm. 1.3.9]{dh}). One then gets that 
  the GIT equivalence induces a stratification of $C^G(X)$ into locally closed  strata (for the euclidean topology on $\NS^G_\bR(X)$; see also \S\ref{semicont} below). To fix the terminology, which is slightly different between \cite{thaddeus} and \cite{dh}, we define:
\begin{definition}\label{defchamber}
A {\it GIT cell} $\sigma$ is a maximal connected locally closed subset of $C^G(X)$ such that  any two (classes of) linearizations $\calL,\calL'\in \sigma\cap \NS^G_\bQ(X)$ are GIT equivalent. {\it A chamber} is a cell of maximal dimension ($=\dim_\bR \NS^G_\bR(X)$). {\it A wall} is a codimension $1$  cell separating two adjacent chambers. 
\end{definition} 

As noted above the stratification is finite and furthermore it is  rational polyhedral in $C^G(X)$. In conclusion, one obtains one of the first major results of VGIT  (cf.  \cite[(2.3), (2.4)]{thaddeus} and \cite[Thm. 0.2.3]{dh}; see also \cite{res}):

\begin{theorem}\label{thmfinite}
 The following hold:
\begin{itemize}
\item[i)] There exist only finitely many cells (and thus also chambers and walls).
\item[ii)] The closure of a cell (in particular, of a chamber) is a convex rational polyhedral cone in $C^G(X)$.
\item[iii)] The cell closures form a fan covering of $C^G(X)$. 
\end{itemize}
\end{theorem}

\begin{remark}\label{dolchamber}
The walls in Dolgachev-Hu \cite{dh} are defined to be the top dimensional strata for which $X^{ss}(\calL)\neq X^{s}(\calL)$.  The chambers of \cite{dh} are the chambers of Def. \ref{defchamber} with the additional condition  $X^{ss}(\calL)=X^{s}(\calL)$. While for many group actions,  the two notions of chambers coincide, there exist examples with $X^{s}(\calL)\subsetneq X^{ss}(\calL)$ for all $\calL\in C^G(X)$ (see  Ressayre's appendix to \cite{dh}) or only for some open regions in the interior of $C^G(X)$ (see \cite{laza}, and \S\ref{sectex} below). 
\end{remark}
 
\subsection{The structure of wall crossings}
\subsubsection{Semi-continuity for the semistable loci and birational transformations.}\label{semicont} A first step towards understanding the change of 
GIT quotient as the linearization moves from one chamber to the other is to notice the following semicontinuity property of the semistable and stable loci (cf.  \cite[Lem. 4.1]{thaddeus}, \cite[\S3.4]{dh}, \cite[Prop. 4]{res}).
\begin{lemma}\label{lemmasemi}
Let $\calL_0\in C^G(X)$ be the class of a $G$-ample linearization. Then for every nearby linearization $\calL\in C^G(X)\cap \NS^G_\bQ(X)\cap B_\epsilon(\calL_0)$ (where $B_\epsilon(\calL_0)$ denotes the Euclidean ball of radius $0<\epsilon\ll 1$ centered at $\calL_0$):
\begin{equation}\label{semieq}
X^s(\calL_0)\subseteq X^{s}(\calL)\subseteq X^{ss}(\calL)\subseteq X^{ss}(\calL_0).
\end{equation}
\end{lemma}

Note that the inclusion of semistable loci from \eqref{semieq} forces the reverse inclusion of the stable loci. Recall that a point is stable w.r.t. $\calL$ iff the stabilizer $G_x$ is finite (an invariant property w.r.t. $\calL$) and $G\cdot x$ is a closed orbit in $X^{ss}(\calL)$. Since $ X^{ss}(\calL_0)$ is an enlargement of $X^{ss}(\calL)$, the orbit $G\cdot x$ might cease to be closed in $X^{ss}(\calL_0)$ due to the inclusion of smaller orbits in the boundary of $G\cdot x$ (see \S\ref{closedorb}). The toy example $\bC\setminus \{0\}=\bC^*\cdot 1\subset \bC=\bC\setminus\{0\}\cup \{0\}$ clearly illustrates this point. In other words, \eqref{semieq} says that in a family of linearizations  $\calL_t$ approaching a special linearization  $\calL_0$ the set of semistable points can jump up, forcing  the quotient corresponding to $\calL_0$ to become ``smaller''. More precisely, in the context of \ref{lemmasemi}, the functorial properties of the GIT quotient give a diagram:
\begin{equation}\label{morphism}
\xymatrix
{ 
X^{s}(\calL)/G\ar@{^{(}->}[d]&X^{s}(\calL_0)/G\ar@{^{(}->}[d]\ar@{_{(}->}[l]\\
X\gquot_\calL G\ar@{->}[r]^{\varphi}&X\gquot_{\calL_0}G}.
\end{equation}
In particular, if $X^s(\calL_0)\neq\emptyset$, the morphism $\varphi$ is clearly birational (with $X^{s}(\calL_0)/G$ a common dense open subset). Furthermore, $\varphi$ is easily described in naive terms. Namely, if $\pi_\calL$ and $\pi_{\calL_0}$ are the two quotient maps and $G\cdot x_0$ is an orbit in $X^{ss}(\calL_0)\setminus X^{ss}(\calL)$, then we have the following contraction
$$\varphi(\pi_\calL(x))=\pi_{\calL_0}(x_0)$$ 
for all $x\in X^{ss}(\calL)$  with closed orbit (and thus separated points in  $X\gquot_\calL G$) such that $\overline{G\cdot x}\supset G\cdot {x_0}$ in $X^{ss}(\calL_0)$. 

\begin{remark}
One should note that $\varphi$ is not always birational, it can be a fibration (but then $\calL_0$ belongs to the boundary of $C^G(X)$). Also, it might happen that $\varphi$ is an isomorphism, even though $X^{ss}(\calL)\neq X^{ss}(\calL_0)$. 
\end{remark}

\subsubsection{The wall crossing transformations in VGIT are often flips.} After the finiteness result of Thm. \ref{thmfinite}, the second set of major results of  VGIT is concerned with structural results on the birational morphism $\varphi$ of \eqref{morphism}. By the general theory of birational transformations, one distinguishes three possibilities for the morphism $\varphi$:
\begin{itemize}
\item[(1)] $\varphi$ is a divisorial contraction;
\item[(2)] $\varphi$ is a small contraction (the exceptional loci have codimension at least $2$);
\item[(3)] $\varphi$ is a fibration.
\end{itemize}
It turns out that the possibilities (1) and (3) are quite special and relatively rare (e.g. (3) only occurs at the boundary of $C^G(X)$). Thus, the interesting case is (2). Again, by   results in birational geometry, it is known that the small contractions force the target space to have bad singularities (e.g. \cite[p. 33]{komo}). The solution to this issue in birational geometry is the {\it flip}, i.e. one replaces a small contraction $f:X\to Y$ by another small contraction $f^+:X^+\to Y$ with certain properties (e.g. \cite[Def. 2.8]{komo}). In particular, if it exists, $X^+$ is uniquely determined and  has good singularities.  

The discussion above suggests that one has to consider not only the morphism $\varphi:X\gquot_\calL G \to X\gquot_{\calL_0} G$ as the stability changes when passing from a chamber to the wall, but also the morphism coming from the other adjacent chamber. In other words, instead of considering the collapsing morphisms as the linearization approaches the  wall, it is better to consider the {\it wall crossing} behavior. More specifically, one considers two (classes of) linearizations $\calL_+, \calL_-\in C^G(X)$ such that the GIT semistability only changes at some $\calL_0\in (\calL_-, \calL_+)$ on the interval joining $\calL_-$ and $\calL_+$ (in $C^G(X)\subset \NS_\bR^G(X)$). The discussion of \S\ref{semicont} automatically gives the following diagram similar to that for flips  from birational geometry:
\begin{equation}\label{flipdiag}
\xymatrix
{ 
X\gquot_- G\ar@{->}[rd]_{\varphi_-}\ar@{-->}[rr]^{g}&&X\gquot_+ G\ar@{->}[ld]^{\varphi_+}\\
&X\gquot_0G&},
\end{equation}
where $\varphi_+, \varphi_-$ are birational morphism as in \eqref{morphism} and $g$ is the induced  birational map (and, for simplicity, we omit $\calL$ from notation). Indeed \eqref{flipdiag} is a flip diagram if one defines a more general notion of flip as follows  (\cite[p. 693]{thaddeus}):
\begin{definition}
Let $f_{-}:X_-\to X$ be a small proper birational morphism. Let $D$ be a $\bQ$-Cartier divisor on $X_{-}$ such that $D$ is relatively negative to $f_{-}$. Then, a $D$-flip is a small proper birational morphism $f_{+}:X_-\to X$ such that 
\begin{itemize}
\item[i)] the Weil divisor $g_*D$ is $\bQ$-Cartier, where $g:X_-\dashrightarrow X_+$ is the induced birational map;
\item[ii)] $g_*D$ is $f_+$-ample. 
\end{itemize}
\end{definition}

We conclude by stating one of the main results of VGIT: the wall crossing birational map $g$ is a flip. More precisely, one has:
\begin{theorem}{\cite[Thm. (3.3)]{thaddeus}}\label{thmflip}
With notation as in \eqref{flipdiag}, if both $X\gquot_- G$ and $X\gquot_+ G$ are not empty, then $\varphi_-$ and $\varphi_+$ are proper and birational. If they are both small, then the rational map $g:X\gquot_- G\dashrightarrow X\gquot_+ G$ is a flip with respect to $\calO(1)$ on $X\gquot_+ G$ (the relative bundle of the projective morphism $X\gquot_+G\to X\gquot_0G$).
\end{theorem}

\begin{remark}
As already hinted, we note that it is possible to have $\varphi_-$ to be a divisorial contraction and $\varphi_+$ to be an isomorphism (see \S\ref{sectex}). However, this is a relatively rare occurrence; most often the wall crossing behavior is flip-like as described in the  theorem.
\end{remark}
 \subsubsection{Explicit description of the flips.}\label{sectexp} The final set of results of VGIT gives a rather explicit description of the morphisms $\varphi_\pm$ and the flip $g$.
 
 We describe first the morphisms $\varphi_{\pm}$ set-theoretically. In the set-up of \eqref{flipdiag}, the relationships \eqref{semieq} among the semistable loci can be made more precise:
 \begin{equation}\label{comparesemiloci}
X^s(\calL_0)=X^s(\calL_-)\cap X^s(\calL_+)\subseteq X^{ss}(\calL_-)\cap X^{ss}(\calL_+)\subsetneq  X^{ss}(\calL_0).
\end{equation}
Let $V=X^{ss}(\calL_+)\cap X^{ss}(\calL_-)$ be the open $G$-invariant subset of $X^{ss}(\calL_0)$ of points that remain semistable when we cross the wall. Clearly, the morphisms $\varphi_{\pm}$ and the rational map $g$ are isomorphisms over the common open subset $U:=V\gquot G$.  The exceptional locus for $\varphi_-$ is then 
$$E_-:=\left(X^{ss}(\calL_-)\setminus X^{ss}(\calL_+)\right)\gquot_{\calL_-} G\subset X\gquot_{\calL_-}G$$
and similarly for $\varphi_+$. It also easy to see that $\varphi_-(E_-)=\varphi_+(E_+)=Z$,
where $Z$ is the closed subset of $X\gquot_0 G$ defined by  
$$Z=\left(X^{ss}(\calL_0)\setminus V\right)\subset  X\gquot_{\calL_0}G.$$
In particular, note that Thm. \ref{thmflip} says that $\dim E_++\dim E_-=\dim X\gquot_0 G-1$. However, one needs to be careful as $Z$ can have several connected components, called {\it strata}, corresponding   to various possibilities for the stabilizers of closed orbits in $X^{ss}(0)\setminus V$.

To understand the scheme theoretic structure of the morphisms $\varphi_{\pm}$, one defines the following ideal sheaves:
$$
\calI_-:=\langle H^0\left(X,\calL_+^N\right)^G\rangle, $$
and similarly $\calI_+$, where $N$ is a sufficiently large and divisible integer. These ideal sheaves descend to ideal sheaves on the quotients $X\gquot_{\pm} G$ which then define the exceptional loci $E_{\pm}$ mentioned above. 
 Similarly, the ideal sheaf corresponding to $Z$ is $(\calI_{+}+\calI_{-})\gquot_0 G$. The following results of Thaddeus then describes the flip $g$ in terms of blow-ups and blow-downs as in the familiar situation of the flip connecting the two small resolution of the three-dimensional quadric cone.

\begin{theorem}[{\cite[Thm. 3.5]{thaddeus}}]\label{thmblow}
The pullbacks of $(\calI_{+}+\calI_{-})\gquot_0 G$ by the morphisms $\varphi_{\pm}:X\gquot_{\pm}G\to X\gquot_0 G$ are exactly $\calI_{\pm}\gquot_{\pm} G$.  The blow-ups of $X\gquot_{\pm}G$ at  $\calI_{\pm}\gquot_\pm G$, and of $X\gquot_0 G$ at $(\calI_{+}+\calI_{-})\gquot_0 G$, are all naturally isomorphic to the irreducible component of the fibered product $X\gquot_{-}G\times_{X\gquot_0 G} X\gquot_{+}G$ dominating $X\gquot_0 G$.
\end{theorem}

Without further simplifying assumptions these blow-ups are quite complicated, for instance they typically involve blow-ups of non-reduced or reducible schemes.  To get more standard blow-ups,  one  makes the following assumptions:
\begin{itemize}
\item[(1)] $X$ is smooth (otherwise the singularities of the quotient will reflect the singularities of $X$); 
\item[(2)] the stabilizers $G_x\cong \Gm$ for points $x\in X^{ss}(\calL_0)\setminus V$ (to guarantee a reasonable structure for the blow-up locus $Z$). 
\end{itemize}
We can now state a final major result of VGIT (cf.  \cite[Thm. 5.6]{thaddeus}, \cite[Thm. 4.2.7]{dh}): 

\begin{theorem}\label{thmexplicit}
Assume that $X$ is smooth and $z\in Z$ corresponds to a closed orbit $G\cdot x$ such that $G_x\cong \Gm$. Then over a neighborhood of $z$ in $Z$, the exceptional divisors $E_{\pm}$ of $\varphi_{\pm}$ are fibrations, locally trivial in the \'etale (or analytic)  topology, with fiber weighted projective spaces.
\end{theorem}

The identification of the weights of the weighted projective fibers in the theorem above can be done using  Luna slice theorem (see \S\ref{sectluna}).

\subsection{Some concluding remarks on VGIT}
We conclude with a summary of the main results of VGIT: 
 \begin{itemize}
 \item[(1)] There are finitely many possibilities for the GIT quotients $X\gquot_\calL G$ as one varies the linearization $\calL$. The set of linearizations is partitioned into rational polyhedral chambers parameterizing GIT equivalent linearizations  (Thm. \ref{thmfinite}).
 \item[(2)]  The semistable loci satisfy a semi-continuity property (Lem. \ref{lemmasemi}). This induces morphisms between quotients for nearby linearizations.
 \item[(3)] The birational change of the GIT quotient as the linearization moves from one chamber to another by passing a wall is a flip (Thm. \ref{thmflip}), which can be understood in terms of blow-ups and blow-downs (Thm. \ref{thmblow}).
 \item[(4)] Under some genericity assumptions, the  flips of (3) can be described quite explicitly (Thm. \ref{thmexplicit}). 
 \end{itemize}
 
These results have  opened the door to a multitude of applications to the construction and geometry of moduli spaces (for a sample see section  \ref{secapp}). More surprisingly, VGIT had strong influences in birational geometry. While the birational geometry seen in VGIT is only a particular case, Hu and Keel \cite{hukeel} have shown that in some sense, it is the most well behaved part of birational geometry. This is discussed in section \ref{sectbir}.

 \subsection{An example of VGIT}\label{sectex}  
 We choose to illustrate the variation of GIT quotients with a somewhat non-standard example (taken from \cite{laza}). We use this example since it exemplifies in a simple geometric situation several aspects of VGIT. Furthermore, it illustrates one of the key strengths of VGIT construction of moduli spaces: its flexibility.  Specifically, one might be interested in two birational models for a moduli space (e.g. obtained by some special constructions). If those two models can be realized as GIT quotients  $X\gquot_{\calL_1}G$ and $X\gquot_{\calL_2}G$ for two choices of linearization,  then, using VGIT, one gets an explicit understanding of the birational relationship between them. For the connection of the VGIT example presented here to other types of constructions of moduli spaces see Prop. \ref{propmmp}, \S\ref{secthodge}, and \cite{laza}.
 
 Our example is probably the simplest VGIT analogue of the plane curve example of \cite[\S4.2]{GIT}. Namely, we consider the moduli space of  pairs $(C,L)$ consisting of a plane curve $C$ of degree $d$ and a line $L$ in $\bP^2$.  The natural GIT set-up for the study of this moduli space is  that of the group $G=\textrm{SL}(3)$ acting diagonally on the parameter space 
  $$X=\bP(H^0(\bP^2,\calO_{\bP^2}(d)))\times\bP(H^0(\bP^2,\calO_{\bP^2}(1)))\cong\bP^N\times \bP^2,$$ 
with $N={d+2 \choose 2}-1$.
In this situation, the space of linearizations $\NS^G(X)$ is identified with $\Pic(X)\cong \bZ\times \bZ$. To fix the notation, let $\pi_1, \pi_2$ be the two projections from $X$ to $\bP^N$ and $\bP^2$ respectively. For $(a,b)\in \bZ\times\bZ$, we define $\calO(a,b)=\pi_1^*\calO(a)\otimes \pi_2^*\calO(b)$ and say  that it has {\it slope} $t=\frac{b}{a}\in \bQ\cup\{\infty\}$. Since replacing a linearization by a multiple does change the GIT quotient, the quotient $X\gquot_\calL G$ and the (semi)stable set $X^{s(s)}(\calL)$ depend only on the slope $t$ of $\calL$; we call a point $x\in X^{s(s)}(\calL)$ {\it $t$-(semi)stable} and denote the quotient by $\calM(t)$.

The nef cone $\overline{A(X)}\subset \NS_\bR^G(X)$ of $X$ is  the upper quadrant $a,b\ge 0$. One then identifies the subcone given by closure of the $G$-ample cone (\cite[(2.5)]{laza}) as:
\begin{equation}
\overline{C^G(X)}=\left\{0\le t=\frac{b}{a}\le \frac{d}{2}, a\ge 0\right\}\subset \overline{A(X)}\subset \NS^G_{\bR}(X)\cong \bR\times\bR.
\end{equation}
The two extremal rays of the cone exemplify the two types of failure of $G$-ampleness. Namely, $t=0$ corresponds to a semi-ample, but not ample linearization (giving the projection $\pi_1:X\to \bP^N$), while $t=\frac{d}{2}$ corresponds to an ample linearization for which there is no stable point. In this case, the GIT quotient still makes sense for these two boundary walls. One then gets:
\begin{itemize}
\item[(1)] $\calM(0)=\bP^N\gquot G$ is the GIT quotient for degree $d$ curves;  the VGIT map $\calM(\epsilon)\to\calM(0)$ is the forgetful map $(C,L)\to C$ (generically a $\bP^2$-fibration).
\item[(2)] $\calM\left(\frac{d}{2}\right)$ is isomorphic to the GIT quotient for $d$ unordered points in $\bP^1$; the VGIT map $\calM\left(\frac{d}{2}-\epsilon\right)\to\calM\left(\frac{d}{2}\right)$ is the forgetful map $(C,L)\to C\cap L\subset L\cong \bP^1$ (generically a weighted projective fibration).
\end{itemize}

This description of the quotient at the boundary walls of $\overline{C^G(X)}$  implies  that the variation of the slope $t\in[0,\frac{d}{2}]$ interpolates between two conditions of stability: the stability
of degree $d$ curves in $\bP^2$ (at $t=0$) and the stability of $d$-tuples of points in $\bP^1$ (at $t=\frac{d}{2}$). In other words, as $t$ increases, $C$ is allowed to be more singular, but the intersection $C\cap L$ should satisfy stronger transversality conditions. This can be made more precise, by relating the GIT approach to the construction of the moduli of pairs to an approach based on MMP as in \cite{hacking}:
\begin{proposition}\label{propmmp}
Let $(C,L)$ be a degree $d$ pair. If the pair $\left(\bP^2,\frac{3}{d+t} (C+t L)\right)$ is log canonical, then $(C,L)$ is $t$-semistable. 
\end{proposition}

\begin{remark}
We emphasize that the implication in the above proposition is only one direction. For some discussion of the relationship between GIT and MMP stability see  \cite[\S10]{hacking} and \cite{kimlee}.
\end{remark}

The finiteness result of VGIT (Thm. \ref{thmfinite}) says that there are only finitely many non-isomorphic GIT quotients $\calM(t)$ and that the cone $\overline{C^G(X)}$ is partitioned into subcones given by the closures of GIT chambers. In this simple situation, there will be a finite number of critical slopes (for which the stability changes) $t_0=0<t_1<\dots<t_k=\frac{d}{2}$ corresponding to the walls. The  subcones partitioning  $\overline{C^G(X)}$ will be spanned  by rays of slopes $t_i$ and $t_{i+1}$.  

For concreteness, we restrict here to pairs of degree $3$. Then, the critical slopes are $t_0=0$, $t_1=\frac{3}{5}$, $t_2=1$,  and $t_3=\frac{3}{2}$.  Furthermore, using the numerical criterion (see \S\ref{sectnum}), one easily computes the stability of degree $3$ pairs: 
\begin{proposition}
Let $(C,L)$ be a degree $3$ pair. If  $L$ passes through a singular point of $C$, then the pair $(C,L)$ is $t$-unstable for all $t>0$. Otherwise, $(C,L)$ is $t$-(semi)stable for
 $t\in (\alpha,\beta)$  (resp. $t\in [\alpha,\beta]$), where $\alpha$ and $\beta$ are given by 
\begin{equation*}
\alpha=\begin{cases}
0           &\textrm{ if }C \textrm{ has at worst nodes}\\
\frac{3}{5} &\textrm{ if }C \textrm{ has an }A_2 \textrm{ singularity}\\
1           &\textrm{ if }C \textrm{ has an }A_3 \textrm{ singularity}\\
\frac{3}{2} &\textrm{ if }C \textrm{ has a }D_4 \textrm{ singularity}                  
\end{cases}
\textrm{ and }
\beta=\begin{cases}
\frac{3}{5} &\textrm{ if }L \textrm{ is inflectional to }C\\
1           &\textrm{ if }L \textrm{ is tangent to }C\\
\frac{3}{2} &\textrm{ if }L \textrm{ is transversal to }C                  
\end{cases}.
\end{equation*}
\end{proposition}

In this example, $X^{ss}(t)=X^s(t)$ for $t\in\left((0,\frac{3}{2})\setminus\{\frac{3}{5},1\}\right)\cap \bQ$, i. e. the quotients corresponding to chambers are geometric quotients.  This is no longer the case starting with degree $4$ (compare Rem. \ref{dolchamber}).  The quotients corresponding to the walls will contain a unique closed orbit of a strictly semistable pair $(C,L)$.  For instance for the wall $t=\frac{3}{5}$, there exists a unique strictly semistable pair $(C_0,L_0)$ with closed orbit which corresponds to $\alpha=\beta=t=\frac{3}{5}$: i.e. a pair consisting of a cuspidal cubic $C_0$ and an inflectional line $L_0$ (in a smooth point of $C_0$). Such a pair is unique up to projective equivalence, and it has a $\Gm$ stabilizer.

As for the birational transformations, note that there exist $7$ quotients: four corresponding to the walls $\calM(0)\cong \bP^1$, $\calM(\frac{3}{5})$, $\calM(1)$, $\calM(\frac{3}{2})\cong \{\mathrm{pt.}\}$ and three geometric quotients corresponding to the chambers $\calM(\epsilon)$, $\calM(1-\epsilon)$, and $\calM(1+\epsilon)$. The boundary morphisms $\calM(\epsilon)\to \calM(0)$ and $\calM(1+\epsilon)\to \calM(\frac{3}{2})$ were already described above in terms of natural forgetful maps. It remains to describe the morphisms at the interior walls.  For $t=1$, one gets $\calM(1-\epsilon)\to \calM(1)$ is a divisorial contraction, and $\calM(1+\epsilon)\cong \calM(1)$. The most interesting case is  $t=\frac{3}{5}$, which gives a flip as in Theorems \ref{thmflip} and \ref{thmexplicit}. The center $Z$ is a point corresponding to the unique closed orbit strictly semistable for $t=\frac{3}{5}$ given by the pair $(C_0,L_0)$ identified in the previous paragraph.  The exceptional divisor $E_+$ parameterizes pairs $(C,L)$ such that  $C$ is a cuspidal cubic and $L$ is a line, not passing through the cusp and not inflectional (see the stability condition from above and \eqref{comparesemiloci}). Similarly,   $E_-$ corresponds to pairs $(C,L)$ such that $C$ is an irreducible (but not cuspidal) cubic and $L$ is an inflectional line. All the orbits of pairs parameterized by $E_+$ and $E_-$ have in their closure the orbit of $(C_0,L_0)$. Thus, the inclusion of $(C_0,L_0)$  in $X^{ss}(\frac{3}{5})$ forces the collapse of $E_{\pm}\cong \bP^1$  to $Z=\{\textrm{pt.}\}$ (via $\varphi_\pm$ as in \eqref{flipdiag}). The local structure of the flip at $t=\frac{3}{5}$ is discussed in \S\ref{sectluna}.

\section{Tools for the analysis of GIT quotients}\label{sectools}
In this section we discuss two essential tools for the study of GIT quotients: {\it the numerical criterion} and {\it Luna's slice theorem}. The numerical criterion gives an efficient way of finding the semistable and stable loci for a group $G$ acting on a projective variety $X$ by reducing the problem to the study of the induced actions for $1$-parameter subgroups of $G$. On the other hand, the slice theorem gives a local description of the quotient by reducing  the action of $G$ to the action of the stabilizer $G_x$ for a closed orbit $G\cdot x\subset X^{ss}$. These tools are well-known and widely used. Here, we focus on  the application of these tools in a VGIT situation. 

\subsection{The Hilbert-Mumford numerical criterion}\label{sectnum}
The numerical criterion is the main tool available for the analysis of the stable and semistable loci in a GIT situation (e.g. see \cite[\S4.2]{GIT} for the simple case of plane curves, and \cite{mumford} for the important case of $\overline{M}_g$). Also, the numerical criterion can be used to prove most of the results of VGIT (\cite{dh}).

Let $\mathcal{L}\in \Pic^G(X)$. For $x\in X$ and $\lambda:\Gm\to G$ a $1$-parameter subgroup ($1$-PS), one defines  $-\mu^\mathcal{L}(x,\lambda)$ to be the weight of the  induced action of $\Gm$ on the fiber $\mathcal{L}_{x_0}\cong \mathbb{A}^1$, where $x_0=\lim_{s\to 0} \lambda(s)\cdot x$. This numerical function $\mu^\mathcal{L}(x,\lambda)$ is used to check the (semi)stability of points in $X$ via the following criterion. 

\begin{theorem}[Hilbert-Mumford Numerical Criterion]\label{numcriterion}
Let $\mathcal{L}$ be an ample $G$-linearized line bundle. Then $x\in X$ is stable 
(resp. semistable) with respect to $\mathcal{L}$ if and only if  $\mu^\mathcal{L}(x,\lambda)>0$ 
(resp. $\mu^\mathcal{L}(x,\lambda)\ge0$) for every nontrivial $1$-PS $\lambda$ of $G$. 
\end{theorem}

The numerical function $\mu^\mathcal{L}(x,\lambda)$ satisfies several properties (e.g. \cite[pg. 49]{GIT}). The most relevant one in the VGIT context is that for fixed $x$ and $\lambda$
$$\mu^\mathcal{L}(x,\lambda):\Pic^G(X)\to \mathbb{Z}$$
is a group homomorphism (N.B. in fact, $\mu^\mathcal{L}$ depends only on the numerical equivalence class of $\calL$ in $\NS^G(X)$). In particular, for two linearizations $\calL_-$ and $\calL_+$, considering the  linearizations on the segment joining them (in $\NS^G_\bQ(X)$)
$$\calL(t)=\calL_{-}^{(1-t)}\otimes \calL_+^{t}, \ \textrm{ for } t\in[0,1]\cap\bQ,$$
one gets
\begin{equation}\label{interpolate}
\mu^{\calL(t)}(x,\lambda)=(1-t)\cdot \mu^{\calL_{-}}(x,\lambda)+t\cdot \mu^{\calL_{+}}(x,\lambda).
\end{equation}
This equation essentially says that the (semi)stability for $\calL(t)$ is an interpolation between the (semi)stability of $\calL_-$ and $\calL_+$ (see \S\ref{sectex} for a geometric example).  Furthermore, the VGIT results (esp. Thm. \ref{thmfinite}) say that the semistability conditions change only at finitely many $t_1,\dots,t_n\in [0,1]\cap\bQ$ (corresponding to intersection of the segement joining $\calL_-$ and $\calL_+$ with walls in $\overline{C^G(X)}$). Thus, a typical analysis of the stability in a VGIT situations involves the following steps:
\begin{itemize}
\item[(1)] describe the (semi)stability for $\calL_{\pm}$;
\item[(2)] identify the critical points $t_1,\dots,t_n$;
\item[(3)] describe the (semi)stability for $\calL(t_i)$ and $\calL(t_i\pm \epsilon)$ (for $0<\epsilon\ll 1$); in particular, identify the blow-up locus $Z(t_i)\subset X\gquot_{\calL(t_i)} G$ and the closed orbits in $X^{ss}(\calL(t_i))$ parameterized by $Z(t_i)$ (see \S\ref{sectexp} for notations).
\end{itemize}
As sketched below, these steps can be accomplished in a somewhat algorithmic way by using the numerical criterion. To understand the flip structure of theorems \ref{thmflip} and \ref{thmexplicit} (and in particular the exceptional loci $E_i^{\pm}$) a more appropriate tool is the slice theorem described in \S\ref{sectluna}. 

\subsubsection{} A GIT analysis for a fixed linearization $\calL$ via the numerical criterion typically has two parts: a combinatorial part and a geometric part.  The combinatorial part consists in fixing a maximal torus $T\subset G$ and  analyzing the stability for the induced $T$ action on $X$ w.r.t. the induced $T$ linearization $\calL_T$. This is what we describe in some detail below. The geometric part consists of extrapolating from $T$-stability to $G$-stability. In practice, this is equivalent to interpreting the $T$-stability in intrinsic geometric terms. More specifically, note 
$$X^{ss}(\calL)=\bigcap_{T' \textrm{ max. torus}} X^{ss}(\calL_{ T'})=\bigcap_{g\in G}X^{ss}(\calL_{g\cdot T\cdot g^{-1}})= \bigcap_{g\in G} g\cdot X^{ss}(\calL_T).$$
Thus, $x$ is $G$-semistable iff all its translates $g\cdot x$ are $T$-semistable. 
Since $x\in X$  represents a geometric object and  $g\cdot x$ is an isomorphic object, one typically aims to show that unstable with respect to $T$ is equivalent to some list of bad geometric features for the object corresponding to $x$. A key idea, due to Kempf \cite{kempf},  is that for an unstable object $x$ there exists an essentially unique maximally destabilizing $1$-PS $\lambda$. This $1$-PS $\lambda$, determines a parabolic subgroup $P_\lambda\subset G$, which is the stabilizer of a flag of linear subspaces (see \cite[\S2.2]{GIT}, \cite{kempf}, \cite[\S9.5]{dolgachev}). The geometric properties that force $x$ to be unstable will be related to this flag (e.g. for plane curves $C$, the destabilizing  flag will typically consists of a singular point $p\in C$ and a special tangent line $L$ through $p$). We emphasize however that, in general, the geometric analysis is quite delicate and it can only be done  case by case. For instance, for hypersurfaces the combinatorial part is quite easy, but the geometric analysis for higher dimensional hypersurfaces was completed only for cubics threefolds and fourfolds (\cite{allcock}, \cite{laza}).

 Now we consider a fixed maximal torus $T\subset G$ and an ample linearization $\calL$ (considered as a $T$-linearization). As described below, the $T$-stability is essentially a combinatorial question. As usual, let $M=\Hom(T,\Gm)\cong \bZ^n$ (where $n$ the dimension of $T$) be the lattice of characters and $N=\Hom(\Gm,T)=M^*$ the dual lattice of $1$-PS.  In particular, there is a natural perfect pairing:
$$\langle\cdot ,\cdot \rangle:M\times N\to \bZ.$$
The linearization $\calL$ gives (after possibly replacing it with a power) a $T$-equivariant embedding $X\hookrightarrow \bP(V)$, with $T$ acting linearly on $V$. In particular, one has an eigenspace decomposition 
\begin{equation}\label{eqdecchi}
V=\bigoplus_{\chi\in M} V_\chi,
\end{equation}
where $V_\chi=\{v\in V\mid t\cdot v=\chi(t)v\}$. Only finitely many characters, say $\chi_1,\dots,\chi_k$ are relevant to the decomposition \eqref{eqdecchi}. An element $x\in X$ has a lift $\tilde x\in V$, which decomposes $\tilde{x}=\sum_{i=1}^k v_{i}$ with $v_{i}\in V_{\chi_i}$. Then the numerical function is simply
$$\mu^\calL(x,\lambda)=\max_{v_i\neq 0} \langle \lambda,\chi_i\rangle.$$
This leads to two basic observations:
\begin{itemize}
\item[i)] for fixed $x$ and $\calL$, $\mu^\calL(x,\lambda)$ is a piecewise linear function in $\lambda\in N\cong \bZ^n$ (here $\lambda$ is $1$-PS varying in  a fixed torus $T$);
\item[ii)] the stability of $x$ actually depends only on the combinatorial object, {\it the state of $x$}: $\st^{\calL}(x)=\{\chi_i\mid v_i\neq 0\}\subset\{\chi_1,\dots,\chi_k\}\subset M\cong \bZ^n$. 
\end{itemize}
To emphasize the second point, we let 
$$\mu(\st^\calL(x),\lambda):=\max_{\chi\in\st^\calL(x) } \langle \chi, \lambda\rangle\left(=\mu^\calL(x,\lambda)\right)$$ 
(see also \cite[p. 306]{kempf}). Since $\mu^\calL(x,\lambda^{-1})=-\min_{v_i\neq 0} \langle \lambda,\chi_i\rangle$, it follows  easily that {\it $x$ is stable is equivalent to the origin in $M_\bR$ being contained in the interior of the convex hull of  $\st^\calL(x)$}. For example, if $X$ is the parameter space for degree $d$ hypersurfaces in $\bP^{n}$, $G=\SL(n+1)$, $T$ is the standard torus, then $V$ is space of degree $d$ polynomials. The  eigenspaces $V_\chi$ are  spanned by the degree $d$ monomials.  The states are then the usual diagrams of Mumford (\cite[p. 10]{mumford}). 

\subsubsection{} Now consider two linearization $\calL^{+}$ and $\calL^{-}$ each assumed very ample. Each of them will give an embedding in a (different) projective space and  an eigenspace decomposition as in \eqref{eqdecchi}. For every $x\in X$, we will get two states $\st^{+}(x)\subset \{\chi^+_1,\dots,\chi^+_{k^+}\}$ and $\st^{-}(x)\subset \{\chi^-_1,\dots,\chi^-_{k^-}\}$ both being  subsets of the character space $M$. For instance, in the situation of \S\ref{sectex}, for a pair $(C,L)$ given by equations $(c,l)$, for appropriate choices, the positive state is the subset of monomials occurring in $c$ and the negative one is the subset of monomials in $l$.

For a line bundle $\calL(t)$ as in \eqref{interpolate}, we  get
 \begin{equation}\label{interpolate2}
\mu^{\calL(t)}(x,\lambda)=(1-t)\cdot \mu(\st^-(x),\lambda)+t\cdot \mu(\st^{+}(x),\lambda).
\end{equation}
We now make the trivial observation that there are only finitely many positive and negative states possible. Denote them by $\st^+_1,\dots, \st^+_{n^+}$ and $\st^-_n,\dots, \st^-_{n^-}$. Furthermore, since $\mu(\st^+_i,\lambda)$   is piecewise linear function in $\lambda$ for each $i\in\{1,\dots,n^+\}$, there exists a finite polyhedral decomposition $\calP$ (with some unbounded regions) of the space of $1$-PS $M$ such that all $\mu(\st^+_i,\lambda)$ and $\mu(\st^-_j,\lambda)$ are simultaneously piecewise linear with respect to $\calP$. From the equation \eqref{interpolate2} the same will be true about the functions $\mu^{\calL(t)}(x,\lambda)$ simultaneously for all $t\in[0,1]$. It follows easily that to test for the positivity of the numerical function $\mu^{\calL(t)}(x,\lambda)$ it suffices to select a finite number of $1$-PS $\lambda_1,\dots,\lambda_s$ that can be determined based on the decomposition $\calP$.

 Since a change of stability at $t\in(0,1)$ corresponds to $\mu^{\calL(t)}(x,\lambda)=0$ and $\mu^{\calL(t-\epsilon)}(x,\lambda)\neq 0$ for some $\lambda$, one gets that all the  critical $t$ (i.e. corresponding to the walls) should be of the form 
$$t_{ij}^{k}=\frac{\mu(\st^{-}_i,\lambda_k)}{\mu(\st^{-}_i,\lambda_k)-\mu(\st^+_j,\lambda_k)}$$ 
with  $\mu(\st^{-}_i,\lambda_k)$ and $\mu( \st^+_j,\lambda_k) $ of opposite signs. The index $i$ runs over the index set $\{1,\dots,n_-\}$ for the negative states. Similarly, $j\in \{1,\dots, n_+\}$, and $k\in\{1,\dots,s\}$. In any case, one obtains a computable finite set of possible critical $t$. In other words, the set $\{\calL(t_{ij}^{k})\}$ contains all the possible walls in the interval that joins $\calL^-$ and $\calL^+$. Not all these values  are actually achieved because of the following two necessary conditions for $\calL(t_{ij}^{k})$ to be realizable as a wall:
\begin{itemize}
\item[a)] the states $\st^{-}_i$ and $\st^{+}_j$ (corresponding to $t_{ij}^k$) should be geometrically realizable for some $x\in X$;
\item[b)]  for $x$ with fixed states $\st^{-}_i$ and $\st^{+}_j$,  one should have $\mu^{\calL(t)}(x,\lambda)\ge 0$ for all $\lambda$ (it suffices to check for $\lambda_i$). 
\end{itemize}

The algorithm outlined above is not very efficient since the number of states $n_{\pm}$ is of order $2^{k_\pm}$, where $k_\pm$ is the number of characters (e.g. degree $d$ monomials) occurring in the decomposition \eqref{eqdecchi} for $\calL_{\pm}$. However, various improvements are well-known (see \cite{bayer}, \cite{morrison}). The author has implemented for \cite{laza} such an improved algorithm that computes all the stability conditions for pairs $(C,L)$ as in \S\ref{sectex} for reasonable degrees $d$ (say $d\le 15$). However, the geometric analysis is feasible only for much lower degrees (say $d\le 6$; compare \cite{shah}).

\subsubsection{} Finally, we note that if $x$ is semistable there exists a $\lambda$ such that 
 $\mu(x,\lambda)=0$ and $x_0=\lim_{s\to 0}\lambda(s)\cdot x$ is still semistable. 
The limit point $x_0$ will be stabilized by a subgroup $\Gm\subset G$ corresponding to the $1$-PS  $\lambda$. Thus, the analysis outlined above  identifies also the relevant closed orbits $G\cdot x_0$ (parametrized by $Z(t)$ at a critical value $t$) as well as their   stabilizers $G_{x_0}$ (or more precisely a maximal torus in $G_{x_0}^0$). 

\subsection{The Luna slice theorem}\label{sectluna}
A good tool for understanding the flips occurring in VGIT (cf. Thm. \ref{thmflip} and Thm. \ref{thmexplicit}) is Luna's slice theorem (\cite{luna}, \cite[Appendix 1.D]{GIT}). 

\begin{notation}
Let $H$ be  a subgroup of $G$. For variety $W$ with a left $H$ action, we denote by $G*_{H} W:=(G\times W)/H$, where
the action of $H$ on the product is given by $h\cdot(g,w)=(gh^{-1},hw)$. Note that  $G*_{H} W$ has a natural left $G$ action. 
\end{notation}

\begin{theorem}[Luna Slice Theorem]\label{lunathm}
Given $x\in X^{ss}$ with closed orbit $G\cdot x$ and $X$ smooth,  there exits a $G_x$-invariant normal slice $V_x\subset X^{ss}$ (smooth and affine) to $G\cdot x$ 
such that we have the following commutative diagram with Cartesian squares: 
$$
\begin{CD}
G*_{G_x}\calN_x@<\text{\'etale}<<G*_{G_x}V_x@>\text{\'etale}>>X^{ss}\\
@VVV                               @VVV@VVV\\
\calN_x/ G_x@<\text{\'etale}<<\left(G*_{G_x}V_x\right)/ G@>\text{\'etale}>> X\gquot G\end{CD}
$$
where $\calN_x$ is the fiber at $x$ of the normal bundle to the orbit $G\cdot x$.
\end{theorem}
\begin{remark}
Note that the stabilizer $G_x$ of a closed orbit in $X^{ss}$ is reductive (Matsushima criterion). Thus, we are still in a standard GIT situation.
\end{remark}

The slice theorem says that a local model for the action of $G$ on $X$ near $x$ is given by the affine quotient $\calN_x\gquot G_x$ (N.B. $\calN_x$ is vector space endowed with the natural $G_x$-action). The theorem can be easily adapted to a VGIT situation. Namely, assume that $G_x^0$ is torus. Then, a local model  for the variation of quotients $X^{ss}\gquot_{\calL_-} G\xrightarrow{\varphi_{-}} X^{ss}\gquot_{\calL_0} G$ around the point $x\in X^{ss}(\calL_0)\setminus X^s(\calL_0)$ with closed orbit is given by $\calN_x\gquot_{\chi_-} {G_x}\to \calN_x/G_x$ from a projective quotient $\calN_x\gquot_{\chi_-} {G_x}$ to the affine quotient $\calN_x/G_x$ for a suitable character $\chi_-$ (compare \S\ref{sectoric}). More specifically, as in the standard case $G_x$ acts on $\calN_x$. Since $G_x$ fixes $x$, $G_x$ acts on the fiber $(\calL_-)_x$ with a character $\chi_-$ (N.B. the character on  $(\calL_0)_x$ is trivial since $x\in X^{ss}(\calL_0)$ and $x$ is stabilized by $G_x$).  In particular, if $G_x\cong \Gm$, the results of Thm. \ref{thmexplicit} are easy to see. Namely, $G_x$ acts with positive, zero, negative weights on the vector space $\calN_x$. With an appropriate choice of $\pm$,  one gets that $\calN_x^{ss}(\chi_-)$ is the complement of the positive weight subspace of $\calN_x$. Thus, the exceptional divisor $E_-$  for $\varphi_{-}$ is a weighted projective bundle over the $0$-weight direction in $\calN_x$ (corresponding to $Z\subset X\gquot_0 G$, see \S\ref{sectexp}) and the weights are the negative weights.

We now  consider the example discussed in \S\ref{sectex}. Namely, we recall that for the wall corresponding to $t=\frac{3}{5}$, there exists a unique closed orbit of a strictly semistable pair $x=(C_0,L_0)$ where $C_0$ is a cubic with a cusp and $L_0$ is inflectional. The defining equations in $\bP^2$ can be taken to be $(x_0x_2^2+x_1^3=0)$ and $(x_0=0)$ respectively. With respect to the chosen coordinates the stabilizer $G_x\cong \Gm$ is diagonal with weights $(5,-1,-4)$. The action of $G_x$ on the normal slice $\calN_x$ is determined  via  the normal bundle sequence:
\begin{equation}\label{normalbundle}
0\xrightarrow{} \calT_{G\cdot x}\xrightarrow{} \calT_{X\mid G\cdot x}\xrightarrow{} \calN_{G\cdot x/X}\xrightarrow{} 0.
\end{equation}
We have $x=(c,l)\in X=\bP^9\times \bP^2=|3L|\times |L|$, which gives the weights of $G_x\cong \Gm$ on $\calT_{X,x}$ are $-18,-12,-9,-6,-3,0,3,6,6,9,9$. Similarly, the weights on $ \calT_{G\cdot x,x}$  are $0,\pm 3, \pm 6, \pm 9$ (N.B. $G=\SL(3)$, $G\cdot x\cong G/G_x$). It follows that the $G_x$ action on $\calN_x$ has weights 
$-18,-12,6,9$.  This gives 
$$E_{-}\cong W\bP(12,18)\cong W\bP(2,3)\cong \bP^1$$ 
and similarly $E_{+}\cong \bP^1$. In conclusion, the birational transformation that occurs at $t=\frac{3}{5}$ for degree $3$ pairs has the effect of replacing the point $Z$ in $X\gquot_{0}G$ corresponding to the pair $(C_0,L_0)$ with the rational curves $E_{\pm}$  in $X\gquot_\pm\calL$ respectively, as described in \S\ref{sectex}.

\section{GIT and birational geometry}\label{sectbir}
A byproduct of VGIT is a treasure of well-behaved examples  in birational geometry, both of local  (flips) and global nature (rational polyhedral decomposition of certain cones). It is perhaps not surprising then that VGIT has  applications (e.g. used in the proof of the weak factorization theorem \cite[\S2.5]{factorization}, \cite{hukeel2}) and influences in birational geometry (e.g. Mori dream spaces). In this section, we review the connections between VGIT and birational geometry.

\subsection{Singularities of GIT quotients} 
Singularities play a central role in birational geometry. We start our discussion of the relationship between GIT and birational geometry with a brief review of the singularities of quotients. First, using the categorical properties of the quotient, it follows that if $X$ is normal, then $X\gquot G$ is also normal. Since the singularities of $X\gquot G$ will reflect the singularities of $X$, we assume for simplicity that $X$ is smooth.  The local structure of the quotient at $\pi(x)$ is described as the quotient $\calN_x/G_x$ of a normal slice (which can be assumed smooth affine) modulo the stabilizer $G_x$ (see \S\ref{sectluna}). Thus, the determining factor for the type of the singularities of the quotient are the stabilizers $G_x$ for $x\in X^{ss}$ with closed orbit. Without any assumptions on stabilizers,   Hochster--Roberts \cite{hochster} proved that $X\gquot G$ has Cohen-Macaulay singularities. Then, Boutot \cite{boutot} strengthened this to $X\gquot G$  has rational singularities (see \cite[Thm. 5.10]{komo} for the relationship to normal Cohen-Macaulay singularities). 

As discussed before, the best  GIT situation is when there are no strictly semistable points ($X^{ss}=X^{s}$). In this case, all the stabilizers are finite and $X\gquot G$ has only finite quotient (or orbifold) singularities. It follows that $X\gquot G$ is  $\bQ$-factorial  with rational singularities (\cite[5.15]{komo}). Furthermore, $X\gquot G$ has log terminal singularities (e.g. \cite[Thm. 2]{schoutens}). Thus, from the point of view of the birational geometry the quotient $X\gquot G$  essentially behaves as a smooth variety.  A more delicate question, which is  relevant in the context of moduli (see \cite{hm} and \cite{hulek}), is whether the singularities of $X\gquot G$ are canonical. This is answered by the Reid-Tai criterion (e.g. \cite[App. 1]{hm}). Namely, the finite quotient singularities are locally of type $\bC^n/G$, for $G$ a finite subgroup of $\GL(n,\bC)$. One can further assume that $G$ acts freely in codimension $1$. Then the singularity is canonical iff for each $g\in G$: $\sum \frac{1}{2\pi i}\log \zeta_k\ge 1$, where $\zeta_k$ are the eigenvalues of $g$, and $\log$ is suitably normalized. In the particular case $G\subset \SL(n,\bC)$ (thus $\prod \zeta_k=1$), the singularities are Gorenstein (i.e. index $1$) canonical singularities (cf.  \cite{watanabe}).

\begin{remark}
For surfaces, the log terminal singularities are precisely the finite quotient singularities. The canonical singularities (i.e. du Val singularities) are those of type $\bC^2/G$ for $G$ a finite subgroup of $\SL(2,\bC)$ (see \cite[(4.18), (4.20)]{komo}).
\end{remark}

Another situation of interest is when $G$ (or every $G_x$) is a torus, then $X\gquot G$ has toric singularities (i.e. the type of singularities that occur for toric varieties). This is a well understood and well behaved class of singularities (see \cite[\S11.4]{coxtoric}). The main issue in this case is that $X\gquot G$ typically fails to be $\bQ$-factorial or even $\bQ$-Gorenstein (if  $X^{ss}\neq X^{s}$). Recall that a toric variety (the local model of the quotient here) is $\bQ$-factorial iff the corresponding fan is simplicial. The more general condition $\bQ$-Gorenstein  has a similar description  (see \cite[11.4.12(a)]{coxtoric}). Thus, one can easily construct examples of quotients $X\gquot G$ that are not $\bQ$-Gorenstein by starting with a fan that does not satisfy the $\bQ$-Gorenstein condition, then the associated toric variety has  a GIT quotient description (see \S\ref{sectoric}) that produces the desired example. On the other hand, assuming that $X\gquot G$ is $\bQ$-Gorenstein, then $X\gquot G$ is automatically log terminal (\cite[Thm. 2]{schoutens}). Furthermore,  assuming $\bQ$-Gorenstein, the conditions of terminal or canonical singularities can be described in terms of the associated fan (\cite[11.4.12(b)]{coxtoric}). 

In conclusion, the singularities of quotients $X\gquot_\calL G$  are mild for linearizations $\calL$ belonging to chambers, but the quotients corresponding to walls typically fail to be $\bQ$-factorial (or even $\bQ$-Gorenstein). Thus,  in order to satisfy the usual assumptions of birational geometry, one  needs to perform a flip (by passing to a nearby chamber). 

\begin{remark} 
For quotients corresponding to linearizations lying on the wall, one can apply Kirwan's (partial) desingularization procedure (\cite{kirwan})  to resolve the singularities corresponding to the closed orbits $G\cdot x$ with $G^0_x$ non-trivial. From the description of the local structure of the quotient (see \S\ref{sectluna}), it is clear that  Kirwan's resolution dominates the VGIT flip  similarly to the case of the cone over the quadric surface (compare with Thm. \ref{thmblow}). 
\end{remark}

\subsection{GIT and toric varieties}\label{sectoric}  
GIT for torus actions on the affine space is essentially equivalent to the theory of toric varieties (\cite{tv}). Consequently, in this situation, GIT and VGIT can be described quite explicitly in combinatorial terms. This leads to numerous non-trivial examples in GIT and birational geometry. Here we only review the basics as a preparation for \S\ref{sectmds}. For more details and examples, we refer the reader to \cite[Ch. 12]{dolgachev}, \cite[Chapters 5, 14,15]{coxtoric}, and \cite[Ch. 6]{mukaib}.

The quotient of an affine space $X=\bA^n$ by the linear action of a torus $G=T\cong(\Gm)^r$ is  an affine toric variety. Namely, it is standard that the action of the torus $T$ on $\bA^n$ can be diagonalized:
$$(t,x)\in T\times \bA^n\to (\chi_1(t)\cdot x_1,\dots,\chi_n(t)\cdot x_n) \in \bA^n,$$
where $\chi_i\in M_T\cong \bZ^r$ are characters of $T$ and $(x_1,\dots,x_n)$ are suitable coordinates on $\bA^n$. Then, $X/G=\Spec R^G$ (see \S\ref{sectaff}) is an affine toric variety. Namely, the ring of invariants is 
$$R^G=k[x_1,\dots,x_n]^T=k[S],$$
where $S$ is the semigroup
$$S=\left\{\sum_{i=1,n} \alpha_i \chi_i=0, \ \alpha_i\in \bZ_{\ge 0} \right\}=\ker\left(M_{\bA^n}\to M_T\right)\cap (\bZ_{\ge 0})^n.$$
Here $M_{\bullet}$ denotes the lattice of characters associated to a toric variety, and the toric structure on $\bA^n$ is determined by the diagonalization of the action of $T$.

More interesting quotients of torus actions on affine spaces are obtained by 
considering projective quotients (see \S\ref{sectproj}).  In this situation, the only line bundle is the trivial bundle $\calO_X$, but there is a choice of linearization  corresponding to a choice $\chi\in M_T$ of  non-trivial character of $T$ (compare Prop. \ref{proplinear}). Similarly to the affine situation, one has $X\gquot_{\chi} G=\Proj(R_\chi^G)$, where the ring $R_\chi^G$  and the grading are defined with respect to $\chi$ as follows:  the degree $d\ge 0$ part of $R_\chi^G$  is $k[S_d]$, where 
$$S_d:=\left\{\sum_{i=1,n} \alpha_i \chi_i=d\chi, \ \alpha_i\in \bZ_{\ge 0} \right\}.$$
Equivalently, $R_\chi^G$ is the usual ring of invariants for the action of $G=T$ on $\bA^{n+1}\cong \bA^n\times \bA^1$ defined by the given action of $T$ on $X=\bA^n$  and by $\chi^{-1}$ on $\bA^1$ (N.B. recall that a projective quotient $X\gquot G$ is defined by considering the affine cone over $X$, here $\bA^{n+1}$; the linearization is a lift of the action on $X$ to the affine cone). It follows easily that $X\gquot_{\chi} G$ is a toric variety, which is projective over the affine toric variety $X/G\cong \Spec k[S_0]$ (compare Rem. \ref{remrelative} and  \cite[Prop. 14.1.12]{coxtoric}). Furthermore, the variation of GIT quotients that occurs when one changes   the linearization on $\calO_{X}$ by means of a character of the torus $T$ is well understood  (in combinatorial terms) and reviewed in other places  (see \cite{ksz},   \cite[Ch. 14, Ch. 15]{coxtoric}). 
 \begin{remark}
 The quotient $X\gquot_\chi G$ is a projective variety if $S_0$ is trivial, which is equivalent to saying that the convex hull of the characters $\chi_i\in M_T\cong \bZ^r$ does not contain the origin (see \cite[Thm. 12.2]{dolgachev}). For example, the weighted projective spaces are GIT quotients of $\bA^{n}$ by a $\Gm$ action with positive weights. 
 \end{remark}

Conversely, any toric variety (projective over affine) can be described as a GIT quotient. In fact, Cox \cite{cox} showed that this can be done in an essentially canonical way. Namely, we recall that a toric variety $V$ (assumed without a torus factor) with associated torus $T_V$  can be described by means of a fan $\Delta$ living in $N_\bR$, where $N=\Hom(\Gm,T_V)$ is the lattice of $1$-parameter subgroups. Let $M=N^*$ be the character lattice. The group of Weil divisors $A_{n-1}(V)$ modulo rational equivalence is described by the exact sequence (\cite[pg. 63]{tv}):
$$0\to M\to \bZ^{\Delta(1)}\to A_{n-1}(V)\to 0$$
where $\bZ^{\Delta(1)}$ is the free abelian group generated by the rays of $\Delta$ that  correspond to $T_V$-invariants Weil divisors. Applying $\Hom(\cdot, \Gm)$ to the above  sequence gives:
\begin{equation}\label{eqtor}
1\to T\to (\Gm)^{\Delta(1)}\to T_V\to 1
\end{equation}
Then the toric variety $V$ has the following presentation as a GIT quotient:
\begin{theorem}[{\cite[Thm. 2.1]{cox}}]\label{thmcox}
Let $V$ be a toric variety determined by a (non-degenerate) fan $\Delta$. Let $X=\bA^{\Delta(1)}$ and $T=\Hom(A_{n-1}(X),\Gm)$ acting on $X$ by means of the inclusion  $T\subset (\Gm)^{\Delta(1)}$. Then $V\cong X\gquot G$. Moreover, the quotient is a geometric quotient (i.e. $X^s=X^{ss}$) iff $V$ is a simplicial toric variety. 
\end{theorem}

\begin{remark} A toric variety has finite abelian quotient singularities iff the corresponding fan $\Delta$ is simplicial. 
\end{remark}

\subsection{Mori dream spaces}\label{sectmds} 
Two of the main results of VGIT  are: the $G$-ample cone $C^G(X)$ is partitioned into finitely many rational polyhedral chambers (Thm. \ref{thmfinite}) and the quotients corresponding to these chambers are related by flips (Thm. \ref{thmflip}). Hu and Keel \cite{hukeel} have noticed that these are intrinsic and very desirable facts about the birational geometry of the quotients $X\gquot G$. 

We recall that cone decompositions occur naturally in birational geometry. Namely, let $V$ be a $\bQ$-factorial projective variety. As usual, $N^1_\bR(V)$ denotes the vector space of numerical equivalence classes  of divisors (with $\bR$ coefficients). Inside $N^1_\bR(V)$, there are two natural cones:
$$\Nef(V)\subseteq \Mov(V)\subset N^1_\bR(V),$$
 {\it the nef cone}  (the closure of {\it the ample cone}) and  (the closure of) {\it the movable cone}. In addition to $\Nef(V)$, the movable cone $\Mov(V)$ contains several other subcones  that are obtained from other birational models $V'$ of $V$. Specifically,  if $f:V\dashrightarrow V'$ is a birational map such that $V'$ is $\bQ$-factorial and $f$ is an isomorphism in codimension $1$, there is a natural identification $f^*:N^1_\bR(V')\cong N^1_\bR(V)$. Then, 
 $f^*(\Nef(V'))$ is a top dimensional subcone of $\Mov(V)$. Furthermore, for non-isomorphic birational models, the associated cones have  
disjoint interiors. Simply note that if $D$ is the pull-back of an ample divisor on $V'$, then  the ring of sections 
$$R(V,D)=\oplus H^0(V,\calO(nD))$$ 
is finitely generated and  $V'\cong \Proj(R(V,D))$. 

Inspired by the VGIT situation, but also for intrinsic reasons, the ideal situation would be that the movable cone decomposes in finitely many rational polyhedral chambers of type $f^*(\Nef(V'))$, called {\it Mori chambers}. Consequently, Hu and Keel \cite{hukeel} have defined the notion of {\it Mori dream space} to capture this situation.

\begin{definition}
A a projective $\bQ$-factorial variety $V$ is  a {\it Mori dream space} if 
\begin{itemize}
\item[i)] $\Pic(V)$ is a finitely generated abelian group (thus $\Pic(V)\otimes \bQ=N^1_\bQ(V)$);
\item[ii)] $\Nef(V)$ is the affine hull of finitely many semi-ample line bundles;
\item[iii)] there are finitely many $f_i:V\dashrightarrow V_i$ which are isomorphisms in codimension one such that each movable divisor $D$ on $X$ is the pullback of some semiample divisor from some model $V_i$. 
\end{itemize}
\end{definition}

As the name suggests, for a Mori dream space $V$ there is a satisfactory understanding  of all birational models of $V$ (see \cite[Prop. 1.11]{hukeel}).
The theory of variations of GIT quotients  produces examples of Mori dream spaces. Namely, generalizing the situation for toric varieties, projective quotients of type $V=X\gquot_\chi T$ (for an affine variety $X$ and  a torus $T$ acting on $X$) are Mori dream spaces (see \cite[Cor. 2.4]{hukeel} for  a precise statement). Varying the choice of linearization  through the character $\chi$, one obtains other birational models for $V$. In fact, under an appropriate identification of the lattice of characters with  $N_\bQ^1(V)$, the Mori chambers in $\Mov(V)$ coincide with the GIT chambers of Thm. \ref{thmfinite} (see \cite[Thm. 2.3]{hukeel}).

More surprisingly, all Mori dream spaces are GIT quotients of the type described above. Again inspired by the case of  toric varieties (see Thm. \ref{thmcox}),  one defines for $V$ a $\bQ$-factorial variety with finitely generated $\Pic(V)$, {\it the Cox ring} 
$$\Cox(V):=\bigoplus H^0(V,L_1^{n_1}\otimes\dots, \otimes L_r^{n_r}),$$
where $L_1,\dots,L_r$ is a basis of $\Pic(V)$.   Then, the main result of Hu and Keel \cite{hukeel} is that Mori dream spaces,  finitely generated Cox rings, and VGIT for torus actions on affine varieties are essentially equivalent notions. 

\begin{theorem}[{\cite{hukeel}}]
Let $V$ be a $\bQ$-factorial projective variety. Then $V$ is a Mori dream space iff the Cox ring $\Cox(V)$ is finitely generated. In this case, $V$ is the quotient of the affine variety $X=\Spec(\Cox (V))$ by the torus $T=\Hom(N^1(V), \Gm)$. The Mori decomposition coincides with the decomposition coming from VGIT by varying the linearization by characters of  $T$.
\end{theorem}
The case of toric varieties discussed in \S\ref{sectoric} corresponds to the case $\Cox(V)$ is the polynomial ring (and thus $X\cong \bA^n$). We refer the reader to \cite[\S15, esp. Thm. 15.1.10]{coxtoric} for a combinatorial description of the nef and moving cones in the case of toric varieties. Other examples of Mori dream spaces include the Fano (and log Fano) varieties  (see \cite[\S1.3]{bchm}).  

For some recent surveys on Mori dream spaces and Cox rings see \cite{moridream} and \cite{laface} respectively. For a more detailed survey on the relationship between birational geometry and VGIT see \cite{hu}.

\section{Applications of birational geometry of GIT quotients to moduli}\label{secapp}
Throughout the history of the subject, GIT and moduli spaces were highly interconnected. As mentioned in the introduction, some of the major achievements of GIT  are the constructions of the  moduli spaces of curves, of abelian varieties, of vector bundles, and of sheaves. Naturally, many of the applications of VGIT have as starting point these successful GIT stories. In particular, a very active and successful area is the application of VGIT to the computation of various cohomological groups associated to certain  moduli spaces of vector bundles.  A good example in this sense (and one of the motivations  of  \cite{thaddeus}) is  Thaddeus' proof of the  Verlinde formula (\cite{thaddeus0}). 
On the other hand, as mentioned in section \ref{sectbir} (esp. \S\ref{sectmds}), some of the applications and influences of VGIT are quite unexpected and not necessarily concerned with moduli spaces. Moreover, the interplay between VGIT and  birational geometry, as well as the recent progress in birational geometry have made the    study of the birational geometry of moduli spaces a central theme of modern moduli theory. 

In this final section, we review two research topics in moduli spaces in which GIT and birational geometry are the central characters.  These two applications are not necessarily the most representative applications of GIT or VGIT, but we believe they illustrate the main points of this survey. Namely, GIT is a useful tool for the construction of moduli spaces and then VGIT  enhances the standard GIT constructions by giving them flexibility. 
 
\subsection{GIT and the birational geometry of $\overline{M}_{g,n}$}
One of the great successes of GIT is its use to prove that the moduli space of curves is a projective variety. We recall that Deligne--Mumford \cite{dm} have constructed a smooth  compactification for the moduli space of curves as a  stack $\overline{\mathcal M}_g$ (a smooth proper Deligne-Mumford stack). Subsequently, Gieseker and Mumford (e.g. \cite{mumford}) have constructed the associated coarse scheme $\overline{M}_g$ of $\overline{\mathcal M}_g$ via GIT by using asymptotic stability for the Hilbert scheme (or Chow variety) of $\nu$-canonical embedded curves ($\nu\gg0$). Subsequently, Mumford and Knudsen \cite{knudsen} have shown the projectivity  of the moduli of pointed stable curves without actually constructing $\overline{M}_{g,n}$ via GIT. This was accomplished only recently by  \cite{baldwin} and \cite{swinarski} (see \cite{msurvey} for a survey on GIT and $\overline{M}_g$). It is worth mentioning that it is possible to construct and prove the projectivity of $\overline{M}_g$ (and $\overline{M}_{g,n}$)  completely avoiding GIT (see for example Koll\'ar \cite{kollarp}). These methods are also applicable to moduli of higher dimensional varieties (e.g. \cite{alexeev}), for which there are very few results obtained via GIT (see however \cite{gieseker} and \cite{viehweg}).  

The birational geometry of $\overline{M}_g$ is a topic of great interest in algebraic geometry ever since the seminal paper of Harris and Mumford \cite{hm}. Initially, the main interest was the Kodaira dimension of $M_g$  and various conjectures on the cones of effective curves and effective divisors on $\overline{M}_g$ (e.g. \cite{hm90}), esp. the Fulton conjecture  (see \cite{gibney}, and \cite{farkas} for a recent survey). Recently, partially inspired by the developments of VGIT, the focus shifted somewhat to the search for various log canonical models for the moduli spaces of curves, pointed curves (see the survey \cite{fs} in this volume), or even stable maps (e.g. \cite{coskun}). Many of the log canonical models of $M_g$  that were constructed so far are obtained via GIT. Here we briefly review two standard examples  $M_{0,n}$ and $M_{g}$ from a GIT perspective (for some results on the birational geometry of $M_{g,1}$ obtained via VGIT, see \cite{jensen}).  

\subsubsection{Birational geometry of $M_{0,n}$} The compactification $\overline{M}_{0,n}$ of the moduli space of ordered $n$ points on $\bP^1$   is easily seen to be birational to $\bP^{n-3}$. Kapranov \cite{kapranov} (and Keel \cite{k1} in a slightly different way) has given an explicit construction for $\overline{M}_{0,n}$ as a sequence of blow-ups  of certain linear configurations in  $\bP^{n-3}$. Moreover, Kapranov gives a description of $\overline{M}_{0,n}$ as the Chow quotient of a Grassmanian by a torus, which can be then related to a VGIT  construction (e.g. \cite[Thm. 1.2]{gmac}). A remaining key question is to understand the cone of effective curves (or equivalently the nef cone) of $\overline{M}_{0,n}$. A conjectural description of the cone of (numerically) effective curves on $\overline{M}_{g,n}$ is given by Fulton's conjecture (see \cite{keelmckernan} included in this volume). Moreover, a positive answer to this conjecture  for $\overline{M}_{0,n}$ would imply the conjecture for all $\overline{M}_{g,n}$ (cf. \cite{gibney}). GIT and VGIT can be  used to understand a (small) slice of the nef cone of $\overline{M}_{0,n}$ (see \cite{alexeevs}). Furthermore, recently discovered connections to conformal blocks and more creative GIT constructions (e.g. \cite{gs}, \cite{giansiracusa}) might shed more light on the Fulton conjecture and related questions. 

More importantly from our perspective, the study of the birational geometry of $\overline{M}_{0,n}$ led to the discovery by Hu and Keel \cite{hukeel} of the close connection between GIT and birational geometry (see \S\ref{sectmds}). Still, the questions asked there (\cite[\S3]{hukeel}): {\it  is $\overline{M}_{0,n}$ a Mori dream space?}, does not seem to have an answer yet (except for $n\le 6$, in which case  $\overline{M}_{0,n}$ is a log Fano variety and thus a Mori dream space; see also \cite{castravet6}).

\subsubsection{Hassett -- Keel program for $M_g$}\label{hkprog}
Another topic of great interest recently is the search for a canonical model for the moduli space of curves $M_g$. We recall that for large $g$ ($g\ge 23$), the moduli of curves $M_g$ is of general type, thus it is a reasonable question to ask for a canonical model $M_g^{can}$. Since $(\overline{M}_g, \Delta)$ (where $\Delta$ denotes the boundary divisor) is a log canonical model, a promising approach is to study the canonical model via the interpolation: 
$$\overline{M}_g(\alpha)=\Proj(\oplus_n H^0(\overline{M}_g, n(K_{\overline{M}_g}+\alpha \Delta)).$$
Note that for $\alpha=1$ one gets $\overline{M}_g$ (for all $g$), while for $g\ge 23$ and $\alpha=0$ one gets $M_g^{can}$. Using general results in birational geometry (esp. \cite{bchm}), one can show that there exists finitely many isomorphism classes $\overline{M}_g(\alpha)$, which are related by birational transformations, giving a picture similar to a VGIT situation. Moreover, it is expected that most of the resulting spaces $\overline{M}_g(\alpha)$ 
 have modular interpretation (giving some alternate compactifications for $M_g$).   The study of the spaces $\overline{M}_g(\alpha)$ and of their modular interpretation is the so called {\it Hassett--Keel program}. Note that even though the ultimate goal of this study is to understand geometrically  $M_g^{can}$, and thus one needs $g\ge 23$, it still makes sense to study $\overline{M}_g(\alpha)$ for small genus.  Results for large $\alpha$ were obtained by Hassett and Hyeon (\cite{hh}, \cite{hh2}) and for low genus by Hassett, Hyeon and Lee (\cite{hg2}, \cite{hl}). A prediction for the critical slopes is given by \cite{alper} (which is closely related to the discussion of \S\ref{sectnum} on finding the critical values in a VGIT situation).

Since the subject is well surveyed in other parts (see esp. \cite{fs} in this volume and \cite{msurvey}), we close by making some brief comments relevant to our survey on applications of GIT to moduli spaces. While Mumford  and Gieseker \cite{mumford} used asymptotic stability on the Hilbert scheme (or Chow variety) of $\nu$-canonically embedded curves for $\nu\ge 5 $ to construct $\overline{M}_g$,  the most powerful tool so far  to construct various $\overline{M}_g(\alpha)$ was to use GIT on Hilbert schemes for $\nu$-canonically embedded curves for small $\nu$. Specifically, if $\Hilb_{g,\nu}^m$ denotes the main component of the scheme of $m$-Hilbert points for $\nu$-canonical curves, then the associated GIT quotients $\Hilb_{g,\nu}^m\gquot \SL(N+1)$ (with $C\xrightarrow{\omega_C^{\otimes\nu}}\bP^N$) tend to produce examples of $\overline{M}_g(\alpha)$. For example, for $\nu=3,4$ and asymptotic linearizations ($m\gg 0$) one gets the moduli of pseudo-stable curves $\overline{M}_g^{ps}\cong \overline{M}_g\left(\frac{9}{11}\right)$ (the first modification of $\overline{M}_g$ as $\alpha$ decreases), see \cite{schubert} and \cite{hm4}.  Similarly,  the case $\nu=2$ and asymptotic linearizations produces the first instance of flip in the Hassett--Keel program (see \cite{hh2}). Finally, only recently it was proved that for $\nu=1,2$ and small $m$ is the generic smooth genus $g$ curve semistable when viewed as a point of $\Hilb_{g,\nu}^m$ (see \cite{afs2}). For a more detailed discussion of the role played by GIT in the Hassett--Keel program see 
\cite[\S2.4]{fs} and \cite{ah}.

 Probably, the two main open questions related to the Hassett--Keel program and GIT are:
\begin{itemize}
\item[(1)] Describe the GIT stability for the Hilbert scheme of canonically embedded curves (see also \cite[\S7.4]{msurvey}).
\item[(2)] A uniform GIT procedure that gives all the spaces $\overline M_g (\alpha)$ as instances of a VGIT problem. 
\end{itemize}
While the behavior of the spaces $\overline M_g (\alpha)$ is as coming from VGIT, the second question seems purely speculative at this point. Also the answer to the first question seems far off with the current techniques. In fact, the questions are interesting even for low genera. Namely, for $g=3$ a canonical (non-hyperelliptic) curve is plane quartic and thus a hypersurface. It follows that there is a unique GIT quotient for canonical genus $3$ curves, which was well understood for a long time (e.g. \cite[p. 80]{GIT}). More recently, a complete analysis of the Hassett--Keel program in this case  was done by Hyeon--Lee \cite{hl}; in particular, the GIT quotient for plane quartics was shown to be the final non-trivial space $\overline{M}_3(\alpha)$. For genus $4$, a canonical curve is a $(2,3)$ complete intersection. A natural parameter space is then a projective bundle $\bP(E)$ over the space of quadrics (N.B. $\bP(E)$ is birational to the corresponding Hilbert scheme $\Hilb_{4,1}$). Since $\Pic(\bP(E))\cong \bZ^2$, GIT for canonical genus $4$ curves leads naturally to a VGIT situation. The resulting VGIT problem and the connection to the Hassett--Keel program are discussed in \cite{g4git} (see also \cite{g4ball} and \cite{fed4}).

\subsection{GIT and Hodge theory}\label{secthodge} 
As mentioned several times in this survey, there are usually several constructions for a moduli space, including GIT. Another standard construction for a moduli space is via Hodge theory (for a discussion of a closely related topic see \cite{milne} in this volume). For instance,   the moduli space of elliptic curves can be described as the quotient $\calh/\SL(2,\bZ)$ of the Siegel upper half space by the modular group.   As mentioned in the introduction there is also a natural GIT construction for the moduli space of elliptic curves. Consequently, one obtains:
\begin{equation}\label{elliptic}
\overline{ M}_1\cong \bP^1\cong \left(\mathfrak h/\SL(2,\bZ)\right)^*\cong \bP\Sym^3 V^*\gquot\SL(3,\bC), 
\end{equation}
where $^*$ denotes the Satake--Baily--Borel compactification (\cite{bb}). This result is somewhat surprising given the different nature of the objects under consideration: $\left(\mathfrak h/\SL(2,\bZ)\right)^*$ is of analytic and arithmetic nature, while  $\bP\Sym^3 V^*\gquot\SL(3,\bC)$ is purely algebraic.

It turns out that \eqref{elliptic} is not an isolated result, but there is a series of similar results: for some of the $M_{0,n}$ with $n\le 12$ (\cite{dmo}), for $M_g$  for $g\le 4$ and $g=6$ (\cite{kondo3}, \cite{l3}, \cite{kondo4}, \cite{ak}), for low degree $K3$ surfaces (\cite{shah}, \cite[\S8.2]{lc2}), and  for moduli of cubic surfaces (\cite{act1}), cubic threefolds (\cite{allcock, act2}, \cite{ls}), and cubic fourfolds (\cite{laza1,laza2} and \cite{looijenga}). More precisely, a similar isomorphism to \eqref{elliptic} only holds for the cases of $M_{0,n}$ considered by Deligne--Mostow \cite{dmo}, and the moduli spaces of cubic surfaces  (\cite{act1}). In all the other cases mentioned above  the GIT construction and the Hodge theoretic construction differ, but  in a rather minimal way: roughly speaking, a Heegner divisor associated to a hyperplane arrangement $\calH$ inside the period domain $\calD/\Gamma$ has to be ``flipped'' (i.e. blown-up to normal crossings and then contracted in the opposite direction). Looijenga \cite{lc1, lc2} has made this statement precise. Namely, for  moduli spaces birational to arithmetic quotients $\calD/\Gamma$ of Type  IV  domains or complex balls, Looijenga gave a comparison theorem (\cite[Thm. 7.6]{lc2}) to GIT quotients, which says that under appropriate hypotheses (satisfied by the geometric examples mentioned above):
\begin{equation}\label{comparison}
\overline{(\calD/\Gamma)}_{\calH}\cong X\gquot_\calL G,
\end{equation}
where $\overline{(\calD/\Gamma)}_{\calH}$ is a birational modification of $\calD/\Gamma$ associated to a hyperplane arrangement $\calH$ (determined by the particular geometric situation). To understand \eqref{comparison}, we recall that the Satake--Baily--Borel compactification $(\calD/\Gamma)^*$ is a projective variety which can be defined as the $\Proj$ of the ring of $\Gamma$-automorphic forms on $\calD$. Similarly, $\overline{(\calD/\Gamma)}_{\calH}$ is $\Proj$ of a ring of meromorphic forms with poles along $\calH$, and thus a projective variety with a tautological polarization. Then, \eqref{comparison} is essentially equivalent to saying that a  linearization $\calL$ and a hyperplane arrangement $\calH$ can be chosen so that the natural polarizations of $X\gquot_\calL G$ and $\overline{(\calD/\Gamma)}_{\calH}$ agree on a open set with high codimension complement.

We note that there are numerous consequences of a result of type \eqref{comparison} for a moduli space. On one hand, the algebraic description as a GIT quotient $X\gquot_\calL G$ can be used to prove properness statements about the period map (e.g. \cite{shah}, \cite{laza2}, \cite{looijenga}). Conversely, the description $\overline{(\calD/\Gamma)}_{\calH}$ comes equipped with a rich arithmetic structure which can be then interpreted geometrically (e.g. results about the Neron--Severi group of $K3$ surfaces).

A particularly interesting case is that of genus $3$ curves. Namely, there exists a natural GIT compactification $\overline{M}^{GIT}_3$ obtained by viewing the smooth non-hyperelliptic curves of genus $3$ as plane quartics, and a  ball quotient description $(\calB_6/\Gamma_6)^*$ due to Kondo \cite{kondo3}. These two birational models of $\overline{M}_3$ are closely related by results of the  type described above. Concretely, there exists a common partial resolution $\widehat{M}_3$ which can be viewed either as a partial Kirwan desingularization of $\overline{M}^{GIT}_3$ or as a Looijenga type arithmetic modification of $(\calB_6/\Gamma_6)^*$ (see  \cite{l3} and  \cite{artebani}). On the other hand, in the context of the Hassett--Keel program (see \S\ref{hkprog}), one studies the log canonical models $\overline{M}_3(\alpha)$ (see \cite{hl}). It turns out that $\overline{M}^{GIT}_3$ and  $(\calB_6/\Gamma_6)^*$ actually occur as the last two non-trivial log canonical models in genus $3$. Specifically, the following holds:
\begin{itemize}
\item[i)] $\overline{M}_3^{GIT}\cong\overline{M}_3(\frac{17}{28})$;
\item[ii)] $\widehat{M}_3\cong \overline{M}_3(\frac{7}{10}-\epsilon)$;
\item[iii)] $(\calB_6/\Gamma_6)^*\cong \overline{M}_3(\frac{7}{10})$.
\end{itemize}
Similar results hold for genus $4$ curves as well (see \cite{kondo4} and \cite{g4ball}).

Another interesting example, where GIT and Hodge theory and also VGIT occur, is  \cite{laza}. Specifically, one considers the moduli space of degree $d$ pairs  $(C,L)$ as described in  \S\ref{sectex}. In the particular case $d=5$, we prove that a special instance of the VGIT quotient (specifically $\calM(1)$ in the notation of \S\ref{sectex}) is isomorphic to the Baily-Borel compactification of an arithmetic quotient $\calD/\Gamma$ of type $IV$. On the other hand, another instance of the quotient (i.e. $\calM(\frac{5}{2}-\epsilon)$) is closely related to the deformation space of a certain class of singularities (namely $N_{16}$). The interest in \cite{laza} is to study this deformation space, but structural results are only known for the space of $\calD/\Gamma$. The VGIT set-up described in \S\ref{sectex} connects these two spaces. In other words, in many situations, VGIT allows one to extract information from a known space (here $\calM(1)\cong(\calD/\Gamma)^*$) and translate it into information about a target space (here $\calM(\frac{5}{2}-\epsilon)$). 
 
\bibliography{vargit}

\newcommand{\etalchar}[1]{$^{#1}$}
\providecommand{\bysame}{\leavevmode\hbox to3em{\hrulefill}\thinspace}
\providecommand{\MR}{\relax\ifhmode\unskip\space\fi MR }
\providecommand{\MRhref}[2]{%
  \href{http://www.ams.org/mathscinet-getitem?mr=#1}{#2}
}
\providecommand{\href}[2]{#2}
\begin{thebibliography}{AKMW02}

\bibitem[ACG{\etalchar{+}}]{abramovich}
D.~Abramovich, Q.~Chen, D.~Gillam, Y.~Huang, M.~Olsson, M.~Satriano, and
  S.~Sun, \emph{{L}ogarithmic {G}eometry and {M}oduli}, to appear in this
  volume (arXiv:1006.5870v1).

\bibitem[ACT]{act2}
D.~Allcock, J.~A. Carlson, and D.~Toledo, \emph{The moduli space of cubic
  threefolds as a ball quotient}, to appear in Memoirs of the A.M.S.

\bibitem[ACT02]{act1}
\bysame, \emph{The complex hyperbolic geometry of the moduli space of cubic
  surfaces}, J. Algebraic Geom. \textbf{11} (2002), no.~4, 659--724.
  \MR{1910264 (2003m:32011)}

\bibitem[AFS10]{alper}
J.~Alper, M.~Fedorchuk, and D.~I. Smyth, \emph{Singularities with
  {$G_m$}-action and the log minimal model program for {$\overline{M}_g$}},
  arXiv:1010.3751v1, 2010.

\bibitem[AFS11]{afs2}
\bysame, \emph{Finite {H}ilbert stability of (bi)canonical curves},
  arXiv:1109.4986, 2011.

\bibitem[AH11]{ah}
J.~Alper and D.~Hyeon, \emph{{GIT} constructions of log canonical models of
  {$M_g$}}, arXiv:1109.2173, 2011.

\bibitem[AK11]{ak}
M.~Artebani and S.~Kond{\=o}, \emph{The moduli of curves of genus six and
  {$K3$} surfaces}, Trans. Amer. Math. Soc. \textbf{363} (2011), no.~3,
  1445--1462.

\bibitem[AKMW02]{factorization}
D.~Abramovich, K.~Karu, K.~Matsuki, and J.~W{\l}odarczyk, \emph{Torification
  and factorization of birational maps}, J. Amer. Math. Soc. \textbf{15}
  (2002), no.~3, 531--572 (electronic). \MR{1896232 (2003c:14016)}

\bibitem[Ale96]{alexeev}
V.~Alexeev, \emph{Moduli spaces {$M_{g,n}(W)$} for surfaces},
  Higher-dimensional complex varieties ({T}rento, 1994), de Gruyter, Berlin,
  1996, pp.~1--22. \MR{1463171 (99b:14010)}

\bibitem[Ale02]{alexeevab}
\bysame, \emph{Complete moduli in the presence of semiabelian group action},
  Ann. of Math. (2) \textbf{155} (2002), no.~3, 611--708. \MR{1923963
  (2003g:14059)}

\bibitem[All03]{allcock}
D.~Allcock, \emph{The moduli space of cubic threefolds}, J. Algebraic Geom.
  \textbf{12} (2003), no.~2, 201--223. \MR{1949641 (2003k:14043)}

\bibitem[Alp08]{alpergit}
J.~Alper, \emph{Good moduli spaces for {Artin} stacks}, arXiv:0804.2242, 2008.

\bibitem[Art09]{artebani}
M.~Artebani, \emph{A compactification of {$M_3$} via {$K3$} surfaces}, Nagoya
  Math. J. \textbf{196} (2009), 1--26. \MR{2591089 (2011a:14070)}

\bibitem[AS08]{alexeevs}
V.~Alexeev and D.~Swinarski, \emph{Nef divisors on {$\overline{M}_{0,n}$} from
  {GIT}}, arXiv:0812.0778, 2008.

\bibitem[AV02]{av}
D.~Abramovich and A.~Vistoli, \emph{Compactifying the space of stable maps}, J.
  Amer. Math. Soc. \textbf{15} (2002), no.~1, 27--75 (electronic). \MR{1862797
  (2002i:14030)}

\bibitem[BB66]{bb}
W.~L. Baily, Jr. and A.~Borel, \emph{Compactification of arithmetic quotients
  of bounded symmetric domains}, Ann. of Math. (2) \textbf{84} (1966),
  442--528. \MR{0216035 (35 \#6870)}

\bibitem[BCHM10]{bchm}
C.~Birkar, P.~Cascini, C.~D. Hacon, and J.~McKernan, \emph{Existence of minimal
  models for varieties of log general type}, J. Amer. Math. Soc. \textbf{23}
  (2010), no.~2, 405--468. \MR{2601039 (2011f:14023)}

\bibitem[BM88]{bayer}
D.~Bayer and I.~Morrison, \emph{Standard bases and geometric invariant theory},
  J. Symbolic Comput. \textbf{6} (1988), no.~2-3, 209--217.

\bibitem[Bou87]{boutot}
J.~F. Boutot, \emph{Singularit\'es rationnelles et quotients par les groupes
  r\'eductifs}, Invent. Math. \textbf{88} (1987), no.~1, 65--68. \MR{877006
  (88a:14005)}

\bibitem[BP90]{brion}
M.~Brion and C.~Procesi, \emph{Action d'un tore dans une vari\'et\'e
  projective}, Operator algebras, unitary representations, enveloping algebras,
  and invariant theory ({P}aris, 1989), Progr. Math., vol.~92, Birkh\"auser
  Boston, Boston, MA, 1990, pp.~509--539. \MR{1103602 (92m:14061)}

\bibitem[BS08]{baldwin}
E.~Baldwin and D.~Swinarski, \emph{A geometric invariant theory construction of
  moduli spaces of stable maps}, Int. Math. Res. Pap. IMRP (2008), no.~1, Art.
  ID rpn 004, 104. \MR{2431236 (2009f:14018)}

\bibitem[Cap94]{caporaso}
L.~Caporaso, \emph{A compactification of the universal {P}icard variety over
  the moduli space of stable curves}, J. Amer. Math. Soc. \textbf{7} (1994),
  no.~3, 589--660. \MR{1254134 (95d:14014)}

\bibitem[Cas09]{castravet6}
A.~M. Castravet, \emph{The {C}ox ring of {$\overline M_{0,6}$}}, Trans. Amer.
  Math. Soc. \textbf{361} (2009), no.~7, 3851--3878. \MR{2491903 (2009m:14037)}

\bibitem[CC]{coskun}
D.~Chen and I.~Coskun, \emph{Towards {M}ori's program for the moduli space of
  stable maps}, to appear in Amer. J. Math. (arXiv:0905.2947).

\bibitem[CLS10]{coxtoric}
D.~A. Cox, J.~Little, and H.~Schenck, \emph{Toric {V}arieties}, to appear in
  Graduate Studies in Mathematics series, 2010.

\bibitem[CMJL12a]{g4ball}
S.~Casalaina-Martin, D.~Jensen, and R.~Laza, \emph{The geometry of the ball
  quotient model of the moduli space of genus four curves}, Compact Moduli
  Spaces and Vector Bundles, Contemp. Math., vol. 564, Amer. Math. Soc.,
  Providence, RI, 2012, pp.~107--136.

\bibitem[CMJL12b]{g4git}
\bysame, \emph{Log canonical models and variation of {GIT} for genus four
  canonical curves}, preprint, 2012.

\bibitem[Cox95]{cox}
D.~A. Cox, \emph{The homogeneous coordinate ring of a toric variety}, J.
  Algebraic Geom. \textbf{4} (1995), no.~1, 17--50. \MR{1299003 (95i:14046)}

\bibitem[DH98]{dh}
I.~V. Dolgachev and Y.~Hu, \emph{Variation of geometric invariant theory
  quotients}, Inst. Hautes \'Etudes Sci. Publ. Math. (1998), no.~87, 5--56.
  \MR{1659282 (2000b:14060)}

\bibitem[DM69]{dm}
P.~Deligne and D.~Mumford, \emph{The irreducibility of the space of curves of
  given genus}, Inst. Hautes \'Etudes Sci. Publ. Math. (1969), no.~36, 75--109.
  \MR{0262240 (41 \#6850)}

\bibitem[DM86]{dmo}
P.~Deligne and G.~D. Mostow, \emph{Monodromy of hypergeometric functions and
  nonlattice integral monodromy}, Inst. Hautes \'Etudes Sci. Publ. Math.
  (1986), no.~63, 5--89. \MR{849651 (88a:22023a)}

\bibitem[Dol03]{dolgachev}
I.~V. Dolgachev, \emph{Lectures on invariant theory}, London Mathematical
  Society Lecture Note Series, vol. 296, Cambridge University Press, Cambridge,
  2003. \MR{2004511 (2004g:14051)}

\bibitem[Far09]{farkas}
G.~Farkas, \emph{Birational aspects of the geometry of {$\overline{ M}_g$}},
  Surveys in differential geometry. {V}ol. {XIV}. {G}eometry of {R}iemann
  surfaces and their moduli spaces, Surv. Differ. Geom., vol.~14, Int. Press,
  Somerville, MA, 2009, pp.~57--110. \MR{2655323}

\bibitem[Fed11]{fed4}
M.~Fedorchuk, \emph{The final log canonical model of the moduli space of stable
  curves of genus four}, arXiv:1106.5012, 2011.

\bibitem[FS]{fs}
M.~Fedorchuk and D.~I. Smyth, \emph{Alternate compactifications of moduli
  spaces of curves}, to appear in this volume (arXiv:1012.0329v1).

\bibitem[Ful93]{tv}
W.~Fulton, \emph{Introduction to toric varieties}, Annals of Mathematics
  Studies, vol. 131, Princeton University Press, Princeton, NJ, 1993, The
  William H. Roever Lectures in Geometry. \MR{1234037 (94g:14028)}

\bibitem[GHS]{hulek}
V.~A. Gritsenko, K.~Hulek, and G.~K. Sankaran, \emph{Moduli of {$K3$}
  {S}urfaces and {I}rreducible {S}ymplectic {M}anifolds}, to appear in this
  volume.

\bibitem[Gia]{giansiracusa}
N.~Giansiracusa, \emph{Conformal blocks and rational normal curves}, to appear
  in J. Algebraic Geom. (arXiv:1012.4835).

\bibitem[Gie77a]{gieseker}
D.~Gieseker, \emph{Global moduli for surfaces of general type}, Invent. Math.
  \textbf{43} (1977), no.~3, 233--282. \MR{0498596 (58 \#16687)}

\bibitem[Gie77b]{giesekervb}
\bysame, \emph{On the moduli of vector bundles on an algebraic surface}, Ann.
  of Math. (2) \textbf{106} (1977), no.~1, 45--60. \MR{466475 (81h:14014)}

\bibitem[GKM02]{gibney}
A.~Gibney, S.~Keel, and I.~Morrison, \emph{Towards the ample cone of
  {$\overline M_{g,n}$}}, J. Amer. Math. Soc. \textbf{15} (2002), no.~2,
  273--294. \MR{1887636 (2003c:14029)}

\bibitem[GM87]{gm}
M.~Goresky and R.~MacPherson, \emph{On the topology of algebraic torus
  actions}, Algebraic groups {U}trecht 1986, Lecture Notes in Math., vol. 1271,
  Springer, Berlin, 1987, pp.~73--90. \MR{911135 (89a:14064)}

\bibitem[GM10]{gmac}
A.~Gibney and D.~Maclagan, \emph{Equations for {C}how and {H}ilbert quotients},
  Algebra Number Theory \textbf{4} (2010), no.~7, 855--885. \MR{2776876}

\bibitem[GS11]{gs}
N.~Giansiracusa and M.~Simpson, \emph{G{IT} compactifications of {$ M_{0,n}$}
  from conics}, Int. Math. Res. Not. IMRN (2011), no.~14, 3315--3334.
  \MR{2817681}

\bibitem[Hac04]{hacking}
P.~Hacking, \emph{Compact moduli of plane curves}, Duke Math. J. \textbf{124}
  (2004), no.~2, 213--257. \MR{2078368 (2005f:14056)}

\bibitem[Has05]{hg2}
B.~Hassett, \emph{Classical and minimal models of the moduli space of curves of
  genus two}, Geometric methods in algebra and number theory, Progr. Math.,
  vol. 235, Birkh\"auser Boston, Boston, MA, 2005, pp.~169--192. \MR{2166084
  (2006g:14047)}

\bibitem[HH08]{hh2}
B.~Hassett and D.~Hyeon, \emph{Log minimal model program for the moduli space
  of curves: The first flip}, arXiv:0806.3444, 2008.

\bibitem[HH09]{hh}
\bysame, \emph{Log canonical models for the moduli space of curves: the first
  divisorial contraction}, Trans. Amer. Math. Soc. \textbf{361} (2009), no.~8,
  4471--4489. \MR{2500894 (2009m:14039)}

\bibitem[HK99]{hukeel2}
Y.~Hu and S.~Keel, \emph{A {GIT} proof of {W{\l}odarczyk} weighted
  factorization theorem}, arXiv:math/9904146, 1999.

\bibitem[HK00]{hukeel}
\bysame, \emph{Mori dream spaces and {GIT}}, Michigan Math. J. \textbf{48}
  (2000), 331--348, Dedicated to William Fulton on the occasion of his 60th
  birthday. \MR{1786494 (2001i:14059)}

\bibitem[HL10]{hl}
D.~Hyeon and Y.~Lee, \emph{Log minimal model program for the moduli space of
  stable curves of genus three}, Math. Res. Lett. \textbf{17} (2010), no.~4,
  625--636. \MR{2661168}

\bibitem[HM82]{hm}
J.~Harris and D.~Mumford, \emph{On the {K}odaira dimension of the moduli space
  of curves}, Invent. Math. \textbf{67} (1982), no.~1, 23--88, With an appendix
  by William Fulton. \MR{664324 (83i:14018)}

\bibitem[HM90]{hm90}
J.~Harris and I.~Morrison, \emph{Slopes of effective divisors on the moduli
  space of stable curves}, Invent. Math. \textbf{99} (1990), no.~2, 321--355.
  \MR{1031904 (91d:14009)}

\bibitem[HM10]{hm4}
D.~Hyeon and I.~Morrison, \emph{Stability of tails and 4-canonical models},
  Math. Res. Lett. \textbf{17} (2010), no.~4, 721--729. \MR{2661175
  (2011f:14077)}

\bibitem[HMSV09]{vakil}
B.~Howard, J.~Millson, A.~Snowden, and R.~Vakil, \emph{The equations for the
  moduli space of {$n$} points on the line}, Duke Math. J. \textbf{146} (2009),
  no.~2, 175--226. \MR{2477759 (2009m:14070)}

\bibitem[HR74]{hochster}
M.~Hochster and J.~L. Roberts, \emph{Rings of invariants of reductive groups
  acting on regular rings are {C}ohen-{M}acaulay}, Advances in Math.
  \textbf{13} (1974), 115--175. \MR{0347810 (50 \#311)}

\bibitem[Hu08]{hu}
Y.~Hu, \emph{Geometric invariant theory and birational geometry}, Third
  {I}nternational {C}ongress of {C}hinese {M}athematicians. {P}art 1, 2, AMS/IP
  Stud. Adv. Math., 42, pt. 1, vol.~2, Amer. Math. Soc., Providence, RI, 2008,
  pp.~155--175. \MR{2409630 (2009i:14063)}

\bibitem[Hum75]{hum}
J.~E. Humphreys, \emph{Linear algebraic groups}, Springer-Verlag, New York,
  1975, Graduate Texts in Mathematics, No. 21. \MR{0396773 (53 \#633)}

\bibitem[Jen]{jensen}
D.~Jensen, \emph{Birational contractions of $\overline{M}_{3,1}$ and
  $\overline{M}_{4,1}$}, to appear in Trans. Am. Math. Soc. (arXiv:1010.3377).

\bibitem[Kap93]{kapranov}
M.~M. Kapranov, \emph{Chow quotients of {G}rassmannians. {I}}, I. {M}.
  {G}el\cprime fand {S}eminar, Adv. Soviet Math., vol.~16, Amer. Math. Soc.,
  Providence, RI, 1993, pp.~29--110. \MR{1237834 (95g:14053)}

\bibitem[Kee92]{k1}
S.~Keel, \emph{Intersection theory of moduli space of stable {$n$}-pointed
  curves of genus zero}, Trans. Amer. Math. Soc. \textbf{330} (1992), no.~2,
  545--574. \MR{1034665 (92f:14003)}

\bibitem[Kem78]{kempf}
G.~R. Kempf, \emph{Instability in invariant theory}, Ann. of Math. (2)
  \textbf{108} (1978), no.~2, 299--316. \MR{506989 (80c:20057)}

\bibitem[Kir84]{kirwancoh}
F.~C. Kirwan, \emph{Cohomology of quotients in symplectic and algebraic
  geometry}, Mathematical Notes, vol.~31, Princeton University Press,
  Princeton, NJ, 1984. \MR{766741 (86i:58050)}

\bibitem[Kir85]{kirwan}
\bysame, \emph{Partial desingularisations of quotients of nonsingular varieties
  and their {B}etti numbers}, Ann. of Math. (2) \textbf{122} (1985), no.~1,
  41--85. \MR{799252 (87a:14010)}

\bibitem[Kir09]{knonred}
\bysame, \emph{Quotients by non-reductive algebraic group actions}, Moduli
  spaces and vector bundles, London Math. Soc. Lecture Note Ser., vol. 359,
  Cambridge Univ. Press, Cambridge, 2009, pp.~311--366. \MR{2537073
  (2011a:14092)}

\bibitem[KL04]{kimlee}
H.~Kim and Y.~Lee, \emph{Log canonical thresholds of semistable plane curves},
  Math. Proc. Cambridge Philos. Soc. \textbf{137} (2004), no.~2, 273--280.
  \MR{MR2090618 (2005m:14055)}

\bibitem[KM]{keelmckernan}
S.~Keel and J.~McKernan, \emph{Contractible extremal rays on
  {$\overline{M}_{0,n}$}}, to appear in this volume (arXiv:alg-geom/9607009v1).

\bibitem[KM76]{knudsen}
F.~F. Knudsen and D.~Mumford, \emph{The projectivity of the moduli space of
  stable curves. {I}. {P}reliminaries on ``det'' and ``{D}iv''}, Math. Scand.
  \textbf{39} (1976), no.~1, 19--55. \MR{0437541 (55 \#10465)}

\bibitem[KM98]{komo}
J.~Koll{\'a}r and S.~Mori, \emph{Birational geometry of algebraic varieties},
  Cambridge Tracts in Mathematics, vol. 134, Cambridge University Press,
  Cambridge, 1998. \MR{1658959 (2000b:14018)}

\bibitem[Kol]{kollar}
J.~Koll\'ar, \emph{Moduli of varieties of general type}, to appear in this
  volume (arXiv:1008.0621v1).

\bibitem[Kol90]{kollarp}
J.~Koll{\'a}r, \emph{Projectivity of complete moduli}, J. Differential Geom.
  \textbf{32} (1990), no.~1, 235--268. \MR{1064874 (92e:14008)}

\bibitem[Kon00]{kondo3}
S.~Kond{\=o}, \emph{A complex hyperbolic structure for the moduli space of
  curves of genus three}, J. Reine Angew. Math. \textbf{525} (2000), 219--232.
  \MR{1780433 (2001j:14039)}

\bibitem[Kon02]{kondo4}
\bysame, \emph{The moduli space of curves of genus 4 and {D}eligne-{M}ostow's
  complex reflection groups}, Algebraic geometry 2000, {A}zumino ({H}otaka),
  Adv. Stud. Pure Math., vol.~36, Math. Soc. Japan, Tokyo, 2002, pp.~383--400.
  \MR{1971521 (2004h:14033)}

\bibitem[KSB88]{ksb}
J.~Koll{\'a}r and N.~I. Shepherd-Barron, \emph{Threefolds and deformations of
  surface singularities}, Invent. Math. \textbf{91} (1988), no.~2, 299--338.
  \MR{922803 (88m:14022)}

\bibitem[KSZ91]{ksz}
M.~M. Kapranov, B.~Sturmfels, and A.~V. Zelevinsky, \emph{Quotients of toric
  varieties}, Math. Ann. \textbf{290} (1991), no.~4, 643--655. \MR{1119943
  (92g:14050)}

\bibitem[Laz09a]{laza}
R.~Laza, \emph{Deformations of singularities and variation of {GIT} quotients},
  Trans. Amer. Math. Soc. \textbf{361} (2009), no.~4, 2109--2161. \MR{2465831
  (2009k:14006)}

\bibitem[Laz09b]{laza1}
\bysame, \emph{The moduli space of cubic fourfolds}, J. Algebraic Geom.
  \textbf{18} (2009), no.~3, 511--545. \MR{2496456 (2010c:14039)}

\bibitem[Laz10]{laza2}
\bysame, \emph{The moduli space of cubic fourfolds via the period map}, Ann. of
  Math. (2) \textbf{172} (2010), no.~1, 673--711. \MR{2680429}

\bibitem[Loo03a]{lc1}
E.~Looijenga, \emph{Compactifications defined by arrangements. {I}. {T}he ball
  quotient case}, Duke Math. J. \textbf{118} (2003), no.~1, 151--187.
  \MR{1978885 (2004i:14042a)}

\bibitem[Loo03b]{lc2}
\bysame, \emph{Compactifications defined by arrangements. {II}. {L}ocally
  symmetric varieties of type {IV}}, Duke Math. J. \textbf{119} (2003), no.~3,
  527--588. \MR{2003125 (2004i:14042b)}

\bibitem[Loo07]{l3}
\bysame, \emph{Invariants of quartic plane curves as automorphic forms},
  Algebraic geometry, Contemp. Math., vol. 422, Amer. Math. Soc., Providence,
  RI, 2007, pp.~107--120. \MR{2296435 (2008b:14045)}

\bibitem[Loo09]{looijenga}
\bysame, \emph{The period map for cubic fourfolds}, Invent. Math. \textbf{177}
  (2009), no.~1, 213--233. \MR{2507640 (2010h:32013)}

\bibitem[LS07]{ls}
E.~Looijenga and R.~Swierstra, \emph{The period map for cubic threefolds},
  Compos. Math. \textbf{143} (2007), no.~4, 1037--1049. \MR{2339838
  (2008f:32015)}

\bibitem[Lun75]{luna}
D.~Luna, \emph{Adh\'erences d'orbite et invariants}, Invent. Math. \textbf{29}
  (1975), no.~3, 231--238. \MR{0376704 (51 \#12879)}

\bibitem[LV09]{laface}
A.~Laface and M.~Velasco, \emph{A survey on {C}ox rings}, Geom. Dedicata
  \textbf{139} (2009), 269--287. \MR{2481851 (2010m:14065)}

\bibitem[Mar77]{maruyama}
M.~Maruyama, \emph{Moduli of stable sheaves. {I}}, J. Math. Kyoto Univ.
  \textbf{17} (1977), no.~1, 91--126. \MR{0450271 (56 \#8567)}

\bibitem[McK10]{moridream}
J.~McKernan, \emph{Mori dream spaces}, Jpn. J. Math. \textbf{5} (2010), no.~1,
  127--151. \MR{2609325}

\bibitem[MFK94]{GIT}
D.~Mumford, J.~Fogarty, and F.~Kirwan, \emph{Geometric invariant theory}, third
  ed., Ergebnisse der Mathematik und ihrer Grenzgebiete (2), vol.~34,
  Springer-Verlag, Berlin, 1994. \MR{1304906 (95m:14012)}

\bibitem[Mil]{milne}
J.S. Milne, \emph{Shimura varieties and moduli}, to appear in this volume
  (arXiv:1105.0887).

\bibitem[Mor09]{msurvey}
I.~Morrison, \emph{G{IT} constructions of moduli spaces of stable curves and
  maps}, Surveys in differential geometry. {V}ol. {XIV}. {G}eometry of
  {R}iemann surfaces and their moduli spaces, Surv. Differ. Geom., vol.~14,
  Int. Press, Somerville, MA, 2009, pp.~315--369. \MR{2655332}

\bibitem[MS11]{morrison}
I.~Morrison and D.~Swinarski, \emph{Gr\"obner techniques for low-degree
  {H}ilbert stability}, Exp. Math. \textbf{20} (2011), no.~1, 34--56.
  \MR{2802723}

\bibitem[Muk03]{mukaib}
S.~Mukai, \emph{An introduction to invariants and moduli}, Cambridge Studies in
  Advanced Mathematics, vol.~81, Cambridge University Press, Cambridge, 2003,
  Translated from the 1998 and 2000 Japanese editions by W. M. Oxbury.
  \MR{2004218 (2004g:14002)}

\bibitem[Muk04]{mukai}
\bysame, \emph{Geometric realization of {$T$}-shaped root systems and
  counterexamples to {H}ilbert's fourteenth problem}, Algebraic transformation
  groups and algebraic varieties, Encyclopaedia Math. Sci., vol. 132, Springer,
  Berlin, 2004, pp.~123--129. \MR{2090672 (2005h:13008)}

\bibitem[Mum77]{mumford}
D.~Mumford, \emph{Stability of projective varieties}, Enseignement Math. (2)
  \textbf{23} (1977), no.~1-2, 39--110. \MR{0450272 (56 \#8568)}

\bibitem[Nag60]{nagata}
M.~Nagata, \emph{On the fourteenth problem of {H}ilbert}, Proc. {I}nternat.
  {C}ongress {M}ath. 1958, Cambridge Univ. Press, New York, 1960, pp.~459--462.
  \MR{0116056 (22 \#6851)}

\bibitem[New78]{newstead}
P.~E. Newstead, \emph{Introduction to moduli problems and orbit spaces}, Tata
  Institute of Fundamental Research Lectures on Mathematics and Physics,
  vol.~51, Tata Institute of Fundamental Research, Bombay, 1978. \MR{546290
  (81k:14002)}

\bibitem[PV89]{popov}
V.~L. Popov and {\`E}.~B. Vinberg, \emph{Invariant theory}, Algebraic geometry,
  4 (Russian), Itogi Nauki i Tekhniki, Akad. Nauk SSSR Vsesoyuz. Inst. Nauchn.
  i Tekhn. Inform., Moscow, 1989, pp.~137--314, 315.

\bibitem[Res00]{res}
N.~Ressayre, \emph{The {GIT}-equivalence for {$G$}-line bundles}, Geom.
  Dedicata \textbf{81} (2000), no.~1-3, 295--324. \MR{1772211 (2001e:14047)}

\bibitem[Sch91]{schubert}
D.~Schubert, \emph{A new compactification of the moduli space of curves},
  Compositio Math. \textbf{78} (1991), no.~3, 297--313. \MR{1106299
  (92d:14018)}

\bibitem[Sch05]{schoutens}
H.~Schoutens, \emph{Log-terminal singularities and vanishing theorems via
  non-standard tight closure}, J. Algebraic Geom. \textbf{14} (2005), no.~2,
  357--390. \MR{2123234 (2006e:13005)}

\bibitem[Ses67]{seshadri}
C.~S. Seshadri, \emph{Space of unitary vector bundles on a compact {R}iemann
  surface}, Ann. of Math. (2) \textbf{85} (1967), 303--336. \MR{0233371 (38
  \#1693)}

\bibitem[Sha80]{shah}
J.~Shah, \emph{A complete moduli space for {$K3$} surfaces of degree {$2$}},
  Ann. of Math. (2) \textbf{112} (1980), no.~3, 485--510. \MR{595204
  (82j:14030)}

\bibitem[Swi]{swinarski}
D.~Swinarski, \emph{{GIT} stability of weighted pointed curves}, to appear in
  Trans. Am. Math. Soc. (arXiv: 0801.1288).

\bibitem[Tha94]{thaddeus0}
M.~Thaddeus, \emph{Stable pairs, linear systems and the {V}erlinde formula},
  Invent. Math. \textbf{117} (1994), no.~2, 317--353. \MR{1273268 (95e:14006)}

\bibitem[Tha96]{thaddeus}
\bysame, \emph{Geometric invariant theory and flips}, J. Amer. Math. Soc.
  \textbf{9} (1996), no.~3, 691--723. \MR{1333296 (96m:14017)}

\bibitem[Vie95]{viehweg}
E.~Viehweg, \emph{Quasi-projective moduli for polarized manifolds}, Ergebnisse
  der Mathematik und ihrer Grenzgebiete (3) [Results in Mathematics and Related
  Areas (3)], vol.~30, Springer-Verlag, Berlin, 1995. \MR{1368632 (97j:14001)}

\bibitem[Wat74]{watanabe}
K.~Watanabe, \emph{Certain invariant subrings are {G}orenstein. {I}, {II}},
  Osaka J. Math. \textbf{11} (1974), 1--8; ibid. 11 (1974), 379--388.
  \MR{0354646 (50 \#7124)}

\end{thebibliography}
\end{document}